\documentclass[3p]{elsarticle}

\usepackage{amsmath}
\usepackage{amsthm}
\usepackage{amssymb}
\usepackage{bm}
\usepackage{accents}
\usepackage{color}
\usepackage{nicematrix}
\usepackage{booktabs}  

\newtheorem{thm}{Theorem}

\newdefinition{rmk}{Remark}
\newproof{pf}{Proof}

\newcommand{\grad}{\mathop{\rm grad}\nolimits}

\numberwithin{equation}{section}

\graphicspath{{figs/}}

\journal{arXiv} 

\begin{document}
	
	\begin{frontmatter}
		
		\title{Computational homogenization of parabolic equations with memory effects for a periodic heterogeneous medium\tnoteref{label1}}
		\tnotetext[label1]{The work was supported by the Russian Science Foundation (grant No. 24-11-00058).}
		
		\author{P.N. Vabishchevich\corref{cor1}\fnref{lab1,lab2}}
		\ead{vab@cs.msu.ru}
		\cortext[cor1]{Corresponding author.}
		
		\address[lab1]{Lomonosov Moscow State University, 1, building 52, Leninskie Gory,  119991 Moscow, Russia}
		
		\address[lab2]{North Caucasus Federal University, 1, Pushkin str., 355017 Stavropol, Russia}
		
		\begin{abstract}
			
			In homogenization theory, mathematical models at the macro level are constructed based on the solution of auxiliary cell problems at the micro level within a single periodicity cell.
			These problems are formulated using asymptotic expansions of the solution with respect to a small parameter, which represents the characteristic size of spatial heterogeneity.
			When studying diffusion equations with contrasting coefficients, special attention is given to nonlocal models with weakly conducting inclusions.
			In this case, macro-level processes are described by integro-differential equations, where the difference kernel is determined by the solution of a nonstationary cell problem.
			The main contribution of this work is the development of a computational framework for the homogenization of nonstationary processes, accounting for memory effects.
			The effective diffusion tensor is computed using a standard numerical procedure based on finite element discretization in space.
			The memory kernel is approximated by a sum of exponentials obtained from solving a partial spectral problem on the periodicity cell.
			The nonlocal macro-level problem is transformed into a local one, where memory effects are incorporated through the solution of auxiliary nonstationary problems.
			Standard two-level time discretization schemes are employed, and unconditional stability of the discrete solutions is proved in appropriate norms.
			Key aspects of the proposed computational homogenization technique are illustrated by solving a two-dimensional model problem.
			
		\end{abstract}
		
		\begin{keyword}
			Nonstationary diffusion equation \sep 
			Homogenization \sep 
			Integro-differential equation \sep
			Exponential sum approximation \sep
			Two-level time discretization schemes
			
			\MSC  35B27 \sep 45K05  \sep 65M12 \sep  65M60   
			
		\end{keyword}
	\end{frontmatter}
	
	\section{Introduction}\label{sec:1}
	
	Modern challenges in mathematical modeling of transport processes in inhomogeneous media require developing effective methods for analysis and numerical solution, particularly for media with pronounced microstructures and complex physical properties. Typical examples include composite materials, porous structures, biological tissues, and geological formations, where transport coefficients may vary significantly in space and time. Direct numerical modeling of such systems faces substantial computational difficulties due to the need to resolve small-scale inhomogeneities, making the application of multiscale methods \cite{efendiev2009multiscale,pavliotis2008multiscale,steinhauser2017computational} particularly relevant.
	
	For periodic structures, fundamental results in homogenization theory have been obtained using two-scale convergence. This concept provides rigorous justification for transitioning from microscopic to effective macroscopic descriptions (see \cite{bakhvalov2012homogenisation,bensoussan2011asymptotic}). The effective diffusion tensor coefficient is computed through solutions of auxiliary cell problems.
	
	Special attention is given to high-contrast diffusion coefficients, where standard homogenization methods require significant modification. For weakly conducting inclusions, memory effects become crucial, and the diffusion process is described by integro-differential equations \cite{sandrakov1,sandrakov2}. The memory kernel is determined by solving nonstationary cell problems on periodic microstructures.
	
	The central challenge addressed in this work involves numerical solutions to integro-differential equations \cite{ChenBook1998}. Conventional approaches requiring complete solution history storage lead to exponential growth in memory and computational resource demands.
	
	More promising computational approaches for memory-effect problems involve transitioning from nonlocal to local formulations \cite{linz1985analytical}. This transformation can be achieved through exponential approximations of memory kernels \cite{vabishchevich2022numerical,vabishchevich2023approximate}. Memory kernel approximation by exponential sums deserves particular attention \cite{braess1986nonlinear}, with fundamental approaches noted in \cite{beylkin2005approximation}.
	
	This paper presents a novel computational homogenization technique for nonstationary processes incorporating memory effects. The technique involves three key steps:
	(i) solving the periodicity cell problem for effective diffusion tensors,
	(ii) approximating memory kernels through spectral problems on cells,
	(iii) constructing extended macroscopic equation systems.
	This approach solves both memory kernel approximation and macroscopic nonlocality problems, showing strong potential for practical applications.
	
	The organization of our work is as follows.
	Section 2 formulates the homogenization problem for nonstationary diffusion equations in periodic media, focusing on memory effects.
	Section 3 describes the nonlocal-to-local transition, determining memory kernels from spectral problems on periodicity cells, enabling direct exponential-sum approximations and extended macroscopic systems.
	Section 4 details the computational algorithm using finite element spatial approximations and two-level temporal schemes, with stability estimates in appropriate spaces.
	Section 5 presents numerical solutions for two-dimensional diffusion problems, demonstrating the algorithm's effectiveness for nonstationary process homogenization.
	Section 6 concludes with study findings.
	
	\begin{figure}[ht]
		\center
		\includegraphics[width=1\linewidth]{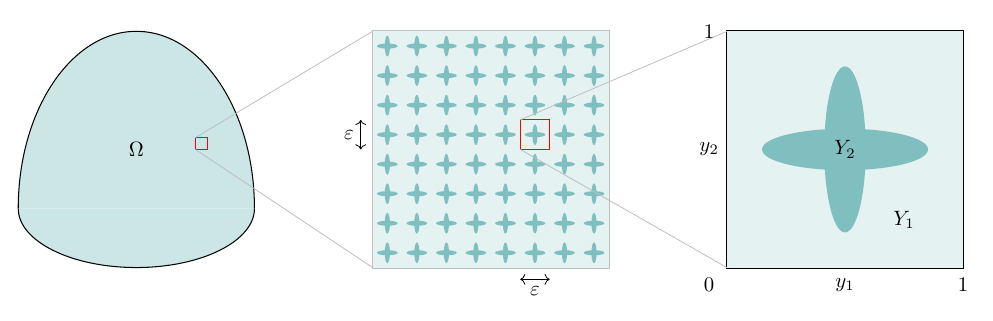}
		\caption{The computational domain $\Omega$ with periodic inhomogeneities of properties and a periodicity cell for a two-component medium.}
		\label{f-1}
	\end{figure}
	
	\section{Problem statement}\label{sec:2}
	
	We consider the problem of nonstationary diffusion in a domain $\Omega \subset \mathbb{R}^2$ (Fig.~\ref{f-1}). The medium is assumed to be periodic with period $\varepsilon$. The substance concentration $u^\varepsilon(x,t)$ is determined from the equation
	\begin{equation}\label{2.1}
		\frac{\partial u^\varepsilon(x,t)}{\partial t} = \nabla \cdot \left( d^\varepsilon(x) \nabla u^\varepsilon(x,t) \right), \quad x \in \Omega,\ t > 0.
	\end{equation}
	The periodically oscillating diffusion coefficient $d^\varepsilon(x)$, which depends on the small parameter $\varepsilon > 0$, has the form:
	\[
	d^\varepsilon(x) = d\left( \frac{x}{\varepsilon} \right),
	\]
	where $d(y)$ is a periodic function of the variable $y = x/\varepsilon$.
	Equation (\ref{2.1}) is supplemented by initial and boundary conditions:
	\begin{equation}\label{2.2}
		u^\varepsilon(x,0) = u_0(x), \quad x \in \Omega, 
	\end{equation}
	\begin{equation}\label{2.3}
		u^\varepsilon(x,t) = 0, \quad x \in \partial\Omega,\ t > 0.
	\end{equation}
	
	We consider two-component media. In this case, the periodicity cell $Y = [0,1]^2$ is divided into two non-overlapping domains: $Y = Y_1 \cup Y_2$ (Fig.~\ref{f-1}).
	The diffusion coefficient $d(y)$ takes constant values in $Y_1$ and $Y_2$. Different homogenization models are used depending on the contrast ratio $\delta$ of the diffusion coefficient in $Y_2$ to that in $Y_1$.
	Analysis of different limiting values of $\delta$ leads to qualitatively different averaged models. Integro-differential homogenization models are associated with the most interesting case of weakly permeable inclusions when $\delta \to 0$.
	We consider the approximation with $\delta \sim \varepsilon^2$:
	\begin{equation}\label{2.4}
		d(y) =
		\begin{cases}
			d_1, & y \in Y_1, \\
			\varepsilon^2 d_2, & y \in Y_2,
		\end{cases}
		\quad \frac{d_1}{d_2} = \mathcal{O}(1).
	\end{equation}
	
	Under the assumptions (\ref{2.4}), the approximate solution to the initial-boundary value problem (\ref{2.1})--(\ref{2.3}) for averaged values was formulated in \cite{sandrakov1,sandrakov2}. The average concentration $u(x,t)$ satisfies the following integro-differential equation:
	\begin{equation}\label{2.5}
		\frac{\partial u(x,t)}{\partial t} + \chi(t) * \frac{\partial u(x,t)}{\partial t} = \nabla \cdot \left( D \nabla u(x,t) \right), \quad x \in \Omega,\ t > 0.
	\end{equation}
	Here $*$ denotes the convolution operator with respect to $t$:
	\[
	\chi(t) * f(t) = \int_0^t \chi(t-s) f(s) \, ds.
	\]
	The boundary and initial conditions for equation (\ref{2.5}) are formulated similarly to (\ref{2.2}), (\ref{2.3}):
	\begin{equation}\label{2.6}
		u(x,0) = u_0(x), \quad x \in \Omega,
	\end{equation}
	\begin{equation}\label{2.7}
		u(x,t) = 0, \quad x \in \partial\Omega,\ t > 0.
	\end{equation}
	The diffusion tensor $D$ and the memory effects kernel $\chi(t)$ in equation (\ref{2.5}) are determined by solving boundary-value problems on the periodicity cell $Y$.
	
	To determine $D$, we use a standard procedure (see, e.g., \cite{hornung1997homogenization}) accounting for the strong contrast in the diffusion coefficient.
	Let $\bm e_i$ be the unit vector in the $i$-th direction, $i=1,2$.
	The auxiliary functions $\theta_i(y)$ are defined as solutions to the equation on the part $Y_1$ of the periodicity cell:
	\begin{equation}\label{2.8}
		\nabla_y \cdot \left( d(y) \left(\bm e_i + \nabla_y \theta_i(y) \right) \right) = 0, \quad y \in Y_1,
		\quad i =1,2.
	\end{equation}
	Periodic boundary conditions are imposed on the external boundary $\partial Y$, while Neumann conditions are specified on the internal boundary:
	\begin{equation}\label{2.9}
		\frac{\partial \theta_i(y)}{\partial n} = - (\bm e_i, \bm n) = 0, \quad y \in \partial Y_2,
		\quad i =1,2,
	\end{equation}
	where $\bm n$ is the outward normal.
	The unique solution is determined by the normalization condition:
	\begin{equation}\label{2.10}
		\int_Y \theta_i(y)\, dy = 0,
		\quad i =1,2.
	\end{equation}
	The effective diffusion tensor $D$ is determined from the solution of the cell problem (\ref{2.8})--(\ref{2.10}):
	\begin{equation}\label{2.11}
		D_{ij} = |Y_1|^{-1}\int_{Y_1} d(y) \left( \delta_{ij} + \frac{\partial \theta_i (y)}{\partial y_j} \right)\, dy,
		\quad i,j = 1,2,
		\quad |Y_1| = \int_{Y_1} dy,
	\end{equation}
	where $\delta_{ij}$ is the Kronecker symbol.
	This tensor is symmetric and positive definite.
	
	For weakly permeable inclusions, we formulate an auxiliary non-stationary problem. 
	The function $\psi(y,t)$ is the solution to the initial-boundary value problem:
	\begin{equation}\label{2.12}
		\frac{\partial \psi(y,t)}{\partial t} = \nabla_y \cdot \left( d_2 \nabla_y \psi(y,t) \right), 
		\quad y \in Y_2,\ t > 0,
	\end{equation}
	\begin{equation}\label{2.13}
		\psi(y,0) = 1, \quad y \in Y_2,
	\end{equation}
	\begin{equation}\label{2.14}
		\psi(y,t) = 0, \quad y \in \partial Y_2,\ t > 0.
	\end{equation}
	The memory kernel $\chi(t)$ is determined from the solution of problem (\ref{2.12})--(\ref{2.14}) by the rule
	\begin{equation}\label{2.15}
		\chi(t) = - |Y_1|^{-1}\int_{Y_2} \frac{\partial \psi(y,t)}{\partial t} \, dy.
	\end{equation}
	
	For numerical solution of the coupled system of equations (\ref{2.5}), (\ref{2.8}), (\ref{2.12}) with (\ref{2.11}), (\ref{2.15}) and the formulated initial and boundary conditions, we can employ standard finite element or finite volume approximations in space \cite{KnabnerAngermann2003,QuarteroniValli1994} and finite difference approximations in time \cite{Ascher2008,SamarskiiTheory}.

	\section{Problem transformation} \label{sec:3}
	
	The main challenges in the approximate solution of the boundary-value problem for equations (\ref{2.5})--(\ref{2.12}) are associated with the integro-differential term in the equation. The difficulties arise from the temporal non-locality, as the solution at the current time depends on its entire history.
	Moreover, the kernel $\chi(t)$ is not given explicitly but is determined according to (\ref{2.15}) from the solution of an auxiliary cell problem. The standard approach \cite{ChenBook1998,mclean1993numerical,mclean1996discretization} based on appropriate quadrature formulas leads to increased computational costs in both time and memory requirements.
	We focus on developing computational algorithms for local mathematical models that account for memory effects. To this end, in our approximate solution we transition from the nonlocal problem by introducing additional unknown quantities \cite{vabishchevich2022numerical,vabishchevich2023approximate}.
	
	\subsection{Auxiliary spectral problem}
	
	For the approximate solution of the problem to determine the kernel (\ref{2.12})--(\ref{2.15}), we consider the spectral problem:
	\begin{equation}\label{3.1}
		\nabla_y \cdot \left( d_2 \nabla_y \varphi(y) \right) + \lambda \varphi(y) = 0, 
		\quad y \in Y_2,
	\end{equation}
	\begin{equation}\label{3.2}
		\varphi(y) = 0, \quad y \in \partial Y_2.
	\end{equation}
	The eigenvalues and eigenfunctions $\lambda_k, \varphi_k(y)$, $k = 1,2, \ldots$ of problem (\ref{3.1})--(\ref{3.2}) are real and satisfy \cite{evans2010partial}:
	\begin{equation}\label{3.3}
		0 < \lambda_1 \leq \lambda_2 \leq \ldots \, .
	\end{equation}
	
	The solution to problem (\ref{2.12})--(\ref{2.14}) is given by
	\begin{equation}\label{3.4}
		\psi(y,t) = \sum_{k=1}^{\infty} (1, \varphi_k)_2 \exp(-\lambda_k t) \varphi_k(y), 
		\quad y \in Y_2,\ t > 0.
	\end{equation}
	Here $(\cdot, \cdot)_2$ denotes the inner product in $L_2(Y_2)$, so that
	\[
	(1, \varphi_k)_2 = \int_{Y_2} \varphi_k(y) dy, 
	\quad k = 1,2, \ldots.
	\]
	
	\subsection{Memory kernel}
	
	Considering (\ref{3.4}), from (\ref{2.15}) we obtain the following representation for the kernel $\chi(t)$:
	\[
	\chi(t) = |Y_1|^{-1} \sum_{k=1}^{\infty} (1, \varphi_k)^2_2 \lambda_k \exp(-\lambda_k t).
	\]
	Thus, the kernel can be expressed as a sum of exponentials:
	\begin{equation}\label{3.5}
		\chi(t) = \sum_{k=1}^{\infty} a_k \exp(-\lambda_k t),
	\end{equation}
	where
	\[
	a_k = |Y_1|^{-1} (1, \varphi_k)^2_2 \lambda_k,
	\quad k = 1,2, \ldots.
	\]
	
	Note that the coefficients in expansion (\ref{3.5}) are positive:
	\begin{equation}\label{3.6}
		a_k > 0, \quad \lambda_k > 0,
		\quad k = 1,2, \ldots.
	\end{equation}
	The kernel $\chi(t)$ is positive-definite \cite{mclean1993numerical,mclean1996discretization} if for $T > 0$ it belongs to $L_1(0, T)$ and satisfies
	\begin{equation}\label{3.7}
		\int_{0}^{T} \zeta(t) \int_{0}^{t} \chi(t-s) \zeta(s) ds dt \geq 0 \quad \forall \zeta \in C[0,T].
	\end{equation}
	The conditions for belonging to $L_1(0, T)$ for (\ref{3.5}) can be verified directly. The sufficient condition for the positive-definite kernel $\chi(t)$ \cite{halanay1965asymptotic}:
	\[
	\chi(t) \geq 0,
	\quad \frac{d\chi}{dt}(t) \leq 0,
	\quad \frac{d^2\chi}{dt^2}(t) \geq 0,
	\quad t > 0,
	\]
	follows from (\ref{3.6}).
	
	Considering this, we establish an a priori estimate for the solution of the homogenized problem (\ref{2.5})--(\ref{2.7}):
	\begin{equation}\label{3.8}
		\|u(x,t)\|_D \leq \|u_0(x)\|_D,
		\quad t > 0,
	\end{equation}
	where
	\[
	\|u(x,t)\|_D^2 = \int_\Omega D \nabla u(x,t) \cdot \nabla u(x,t) \, dx.
	\]
	We will use the stability estimate with respect to initial data (\ref{3.8}) as a basis for considering time discretizations. For the $L_2(\Omega)$ norm, we use the notation $\|\cdot\|$.
	
	\subsection{Approximate solution of the homogenized problem}
	
	Based on the kernel representation as a sum of exponentials (\ref{3.5}), we can construct a computational algorithm for solving our nonlocal problem. We focus on methods that approximate the kernel by a finite sum, keeping only the first $m$ eigenvalues. 
	In representation \eqref{3.5}, we set
	\[
	\chi(t) = \chi_m(t) + \widetilde{\chi}_m(t),
	\]
	where
	\[
	\begin{split}
		\chi_m(t) & = \sum_{k=1}^{m} a_k \exp(-\lambda_k t), \\
		\widetilde{\chi}_m(t) & = |Y_1|^{-1} \sum_{k=m+1}^{\infty} (1, \varphi_k)^2_2 \lambda_k \exp(-\lambda_k t).
	\end{split}
	\]
	Special attention should be paid to estimating the influence of the term $\widetilde{\chi}_m(t)$. 
	
	For sufficiently large $m$, the eigenvalues $\lambda_k$, $k = m+1, m+2, \ldots$ are large. The function
	\[
	s(t) = 
	\begin{cases}
		\lambda \exp(-\lambda t), & t \geq 0, \\
		0, & t < 0
	\end{cases}
	\]
	approaches the delta function as $\lambda \rightarrow \infty$. Considering this, we obtain
	\[
	\widetilde{\chi}_m(t) \approx r_m \delta(t),
	\quad r_m = |Y_1|^{-1} \sum_{k=m+1}^{\infty} (1, \varphi_k)^2_2.
	\]
	Using the introduced notation, we have
	\begin{equation}\label{3.9}
		r_m = (1-|Y_2|)^{-1} \left(|Y_2| - \sum_{k=1}^{m} (1, \varphi_k)^2_2 \right).
	\end{equation}
	The corresponding approximate solution $\widetilde{u}(x,t)$ is obtained as the solution to the boundary-value problem:
	\begin{equation}\label{3.10}
		(1 + r_m) \frac{\partial \widetilde{u}(x,t)}{\partial t} + \chi_m(t) * \frac{\partial \widetilde{u}(x,t)}{\partial t} = \nabla \cdot \left( D \nabla \widetilde{u}(x,t) \right), \quad x \in \Omega,\ t > 0,
	\end{equation}
	\begin{equation}\label{3.11}
		\widetilde{u}(x,0) = u_0(x), \quad x \in \Omega,
	\end{equation}
	\begin{equation}\label{3.12}
		\widetilde{u}(x,t) = 0, \quad x \in \partial\Omega,\ t > 0.
	\end{equation}
	
	To formulate a local problem, we introduce auxiliary functions:
	\begin{equation}\label{3.13}
		v_k(x,t) = \int_0^t \exp(-\lambda_k (t-s)) \frac{\partial \widetilde{u}}{\partial s}(x,s) ds,
		\quad k=1,2,\ldots,m.
	\end{equation}
	In terms of these functions, equation (\ref{3.10}) can be rewritten as
	\begin{equation}\label{3.14}
		(1 + r_m) \frac{\partial \widetilde{u}(x,t)}{\partial t} + \sum_{k=1}^m a_k v_k(x,t) = \nabla \cdot \left( D \nabla \widetilde{u}(x,t) \right), \quad x \in \Omega,\ t > 0.
	\end{equation}
	The introduced auxiliary functions $v_k(x,t)$, $k = 1,2,\ldots,m$, satisfy the equations:
	\begin{equation}\label{3.15}
		\frac{\partial v_k(x,t)}{\partial t} + \lambda_k v_k(x,t) - \frac{\partial \widetilde{u}(x,t)}{\partial t} = 0,
		\quad k = 1,2,\ldots,m.
	\end{equation}
	The system of equations (\ref{3.15}) is supplemented by the initial conditions:
	\begin{equation}\label{3.16}
		v_k(x,0) = 0,
		\quad x \in \Omega,
		\quad k = 1,2,\ldots,m.
	\end{equation}
	We rewrite the equation in a more convenient form. From (\ref{3.15}) we obtain
	\[
	v_k(x,t) = \frac{1}{\lambda_k} \frac{\partial \widetilde{u}(x,t)}{\partial t} - \frac{1}{\lambda_k} \frac{\partial v_k(x,t)}{\partial t},
	\quad k = 1,2,\ldots,m.
	\]
	Substituting into (\ref{3.14}) yields
	\begin{equation}\label{3.17}
		\left(1 + r_m + \sum_{k=1}^{m} \frac{a_k}{\lambda_k}\right) \frac{\partial \widetilde{u}(x,t)}{\partial t}
		- \sum_{k=1}^{m} \frac{a_k}{\lambda_k} \frac{\partial v_k(x,t)}{\partial t} = \nabla \cdot \left( D \nabla \widetilde{u}(x,t) \right), \quad x \in \Omega,\ t > 0.
	\end{equation}
	
	\begin{thm}\label{t-1}
		For the solution of problem (\ref{3.11}), (\ref{3.12}), (\ref{3.5})--(\ref{3.17}), the following a priori stability estimate with respect to initial data holds:
		\begin{equation}\label{3.18}
			\|\widetilde{u}(x,t)\|^2_D + \sum_{k=1}^{m} a_k \|v_k(x,t)\|^2 \leq \|u_0(x)\|^2_D,
			\quad t > 0.
		\end{equation}
	\end{thm}
	
	\begin{pf}
		We multiply equation (\ref{3.17}) by $\partial \widetilde{u}(x,t)/\partial t$ and integrate over $\Omega$, while multiplying each equation (\ref{3.15}) for $k = 1,2,\ldots,m$ by $a_k \lambda_k^{-1} \partial v_k(x,t)/\partial t$. This yields the equalities:
		\[
		(1 + r_m) \left\|\frac{\partial \widetilde{u}(x,t)}{\partial t}\right\|^2 + 
		\frac{1}{2} \frac{d}{dt} \|\widetilde{u}(x,t)\|^2_D + 
		\sum_{k=1}^{m} \frac{a_k}{\lambda_k} \left(\frac{\partial \widetilde{u}(x,t)}{\partial t}, \frac{\partial \widetilde{u}(x,t)}{\partial t}\right) -
		\sum_{k=1}^{m} \frac{a_k}{\lambda_k} \left(\frac{\partial v_k(x,t)}{\partial t}, \frac{\partial \widetilde{u}(x,t)}{\partial t}\right) = 0,
		\]
		\[
		\frac{a_k}{\lambda_k} \left(\frac{\partial v_k(x,t)}{\partial t}, \frac{\partial v_k(x,t)}{\partial t}\right) -
		\frac{a_k}{\lambda_k} \left(\frac{\partial \widetilde{u}(x,t)}{\partial t}, \frac{\partial v_k(x,t)}{\partial t}\right) +
		a_k \frac{1}{2} \frac{d}{dt} \|v_k(x,t)\|^2 = 0,
		\quad k = 1,2,\ldots,m.
		\]
		Adding them, we obtain
		\[
		(1 + r_m) \left\|\frac{\partial \widetilde{u}(x,t)}{\partial t}\right\|^2 + 
		\frac{1}{2} \frac{d}{dt} \left(\|\widetilde{u}(x,t)\|^2_D + 
		\sum_{k=1}^{m} {a_k} \|v_k(x,t)\|^2\right) + 
		\sum_{k=1}^{m} \frac{a_k}{\lambda_k} \left\| \frac{\partial v_k(x,t)}{\partial t} - \frac{\partial \widetilde{u}(x,t)}{\partial t} \right\|^2 = 0.
		\]    
		This leads us to the inequality
		\[
		\frac{d}{dt} \left(\|\widetilde{u}(x,t)\|^2_D + \sum_{k=1}^{m} a_k \|v_k(x,t)\|^2\right) \leq 0,
		\]
		from which the desired estimate (\ref{3.18}) follows.
		\qed
	\end{pf}
	
	\section{Computational algorithm} \label{sec:4}
	
	The proposed computational homogenization is based on sequentially solving the following problems:
	
	\begin{itemize}
		\item Numerical solution of boundary-value problems \eqref{2.8}--\eqref{2.10} on the periodicity cell $Y_1$ to determine the effective diffusion tensor according to \eqref{2.11};
		\item Calculation of coefficients in the approximation \eqref{3.9} for the memory kernel based on solving the partial eigenvalue problem \eqref{3.1}, \eqref{3.2} $(\lambda_k, \varphi_k(y)$, $k = 1,2,\ldots,m)$ in domain $Y_2$;
		\item Solution of problem (\ref{3.11}), (\ref{3.12}), (\ref{3.5})--(\ref{3.17}) to find the homogenized solution in the computational domain $\Omega$.
	\end{itemize}
	We highlight the key elements of these problems.
	
	\subsection{Effective diffusion coefficient}
	
	The numerical solution of the problems is based on using piecewise linear finite elements for spatial approximation on triangular grids. The approximate solution of problems \eqref{2.8}--\eqref{2.10} in domain $Y_1$ is denoted by $\theta_i^h(y)$, $i = 1,2$. 
	On the set of periodic functions $V_1^h \subset H^{1}(Y_1)$, we seek $\theta_i^h(y) \in V_1^h$, $i = 1,2$ satisfying
	\begin{equation} \label{4.1}
		(d \grad_y \theta_i^h \cdot \grad_y \xi_1)_1 + \mu (1, \xi_1)_1 + \int_{\partial Y_2} d \frac{\partial \theta_i^h}{\partial y_i} \xi_1 dy = 0
		\quad \forall w_1 \in V_1^h, \quad i = 1,2.
	\end{equation}
	Here $(\cdot, \cdot)_1$ denotes the scalar product in $L_2(Y_1)$, and $\mu$ is the Lagrange multiplier for enforcing condition \eqref{2.10}.
	After solving problem \eqref{4.1}, the components of the diffusion tensor $D^h$ are calculated according to \eqref{2.11}:
	\begin{equation}\label{4.2}
		D_{ij}^h = |Y_1|^{-1}\int_{Y_1} d(y) \left( \delta_{ij} + \frac{\partial \theta_i^h (y)}{\partial y_j} \right) dy,
		\quad i,j = 1,2.
	\end{equation}
	Problem \eqref{4.1}, \eqref{4.2} is solved on the part $Y_1$ of the periodicity cell $Y$.
	
	\subsection{Memory kernel calculation}
	
	First, a partial spectral problem is solved in the part $Y_2$ of the periodicity cell. 
	We seek an approximate solution $\varphi^h(y)$ from the class of piecewise linear functions $V_2^h \subset H^{1}(Y_2)$ that satisfy homogeneous Dirichlet boundary conditions on $\partial Y_2$. Considering \eqref{3.1}, \eqref{3.2} and the introduced notation, we have
	\begin{equation} \label{4.3}
		(d_2 \nabla_y \varphi^h, \nabla_y \xi_2)_2 + \lambda^h (\varphi^h, \xi_2)_2 = 0 
		\quad \forall \xi_2 \in V_2^h.
	\end{equation}
	Numerically, we determine the first $m$ eigenvalues and eigenfunctions $\lambda_k^h, \varphi_k^h(y)$, $k = 1,2,\ldots,m$ of problem (\ref{4.3}).
	
	According to \eqref{3.5}, the memory kernel is approximated by a sum of exponentials:
	\begin{equation}\label{4.4}
		\chi_m^h(t) = \sum_{k=1}^{m} a_k^h \exp(-\lambda_k^h t),
	\end{equation}
	where
	\[
	a_k^h = (1-|Y_2|)^{-1} (1, \varphi_k^h)^2_2 \lambda_k^h,
	\quad k = 1,2,\ldots,m.
	\]
	The approximation for the term $r_m$ according to \eqref{3.9} is given by
	\[
	r^h_m = (1-|Y_2|)^{-1} \left(|Y_2| - \sum_{k=1}^{m} (1, \varphi^h_k)^2_2 \right).
	\]
	For problem \eqref{4.3}, \eqref{4.4}, we employ well-developed computational algorithms \cite{saad2011numerical} for finding the first minimal (in absolute value) eigenvalues of generalized problems for symmetric matrices.
	
	\subsection{Numerical solution of the homogenized problem}
	
	A triangular mesh is introduced in the computational domain $\Omega$. For approximation by linear finite elements on this mesh, we define the space of functions $V^h \subset H^1(\Omega)$ that satisfy homogeneous (see \eqref{3.12}) Dirichlet boundary conditions. We seek an approximate solution to problem (\ref{3.11}), (\ref{3.12}), (\ref{3.5})--(\ref{3.17}) $y(x,t) \in V^h, w_k(x,t) \in V^h, \ k = 1,2,\ldots,m$, corresponding to $\widetilde{u}(x,t)$ and $v_k(x,t), \ k = 1,2,\ldots,m$.
	We need to find $y(x,t) \in V^h$ and $w_k(x,t) \in V^h, \ k = 1,2,\ldots,m$ satisfying:
	\begin{equation}\label{4.5}
		\left(1 + r^h_m + \sum_{k=1}^{m} \frac{a_k}{\lambda_k}\right)
		\left(\frac{\partial y}{\partial t}, z\right) - 
		\sum_{k=1}^{m} \frac{a_k}{\lambda_k} \left(\frac{\partial w_k}{\partial t}, z\right) + \left(D \nabla y, \nabla z\right) = 0,
	\end{equation} 
	\begin{equation}\label{4.6}
		\left(\frac{\partial w_k}{\partial t}, z_k\right) + \lambda_k (w_k, z_k) -
		\left(\frac{\partial y}{\partial t}, z_k\right) = 0,
		\quad k = 1,2,\ldots,m, 
		\quad t > 0,
	\end{equation} 
	\begin{equation}\label{4.7}
		(y(x,0), z(x)) = (u_0(x), z(x)),
	\end{equation}
	\begin{equation}\label{4.8}
		(w_k(x,0), z_k(x)) = 0,
		\quad k = 1,2,\ldots,m,
	\end{equation}
	where $(\cdot, \cdot)$ denotes the scalar product in $L_2(\Omega)$, 
	for all $z(x), z_k(x) \in V^h, \ k = 1,2,\ldots,m$.
	
	By choosing in \eqref{4.5}, \eqref{4.6}
	\[
	z = \frac{\partial y}{\partial t},
	\quad z_k = \frac{a_k}{\lambda_k} \frac{\partial w_k}{\partial t},
	\quad k = 1,2,\ldots,m,
	\]
	and following the proof of Theorem~\ref{t-1}, we establish the estimate \eqref{3.18} for the solution of problem \eqref{4.5}--\eqref{4.8}.
	Special attention should be paid to the time approximation problem.
	
	We use a uniform temporal grid with step size $\tau$: $t_n = n\tau$, $n = 0,1,\ldots$. The approximate solution at time $t_n$ is denoted by $y^n$, $w_k^n$, $k = 1,2,\ldots,m$.
	For the approximate solution of problem \eqref{4.5}--\eqref{4.8}, we consider standard two-level approximations with a weighting parameter $\sigma$ \cite{SamarskiiTheory,SamarskiiMatusVabischevich2002}:
	\begin{equation}\label{4.9}
		\left(1 + r^h_m + \sum_{k=1}^{m} \frac{a_k}{\lambda_k}\right)
		\left(\frac{y^{n+1} - y^{n}}{\tau}, z\right) - 
		\sum_{k=1}^{m} \frac{a_k}{\lambda_k} \left(\frac{w_k^{n+1}-w_k^{n}}{\tau}, z\right) + \left(D \nabla y^{n+\sigma}, \nabla z\right) = 0,
	\end{equation} 
	\begin{equation}\label{4.10}
		\left(\frac{w_k^{n+1} - w_k^{n}}{\tau}, z_k\right) + \lambda_k (w_k^{n+\sigma}, z_k) -
		\left(\frac{y^{n+1} - y^{n}}{\tau}, z_k\right) = 0,
		\quad k = 1,2,\ldots,m, 
		\quad n = 0,1,\ldots,
	\end{equation} 
	\begin{equation}\label{4.11}
		(y^0, z) = (u_0, z),
	\end{equation}
	\begin{equation}\label{4.12}
		(w_k^0, z_k) = 0,
		\quad k = 1,2,\ldots,m,
	\end{equation}
	using the notation
	\[
	y^{n+\sigma} = \sigma y^{n+1} + (1-\sigma) y^{n}.
	\]
	For $\sigma = 1/2$, the difference scheme \eqref{4.9}--\eqref{4.12} has second-order approximation in $\tau$, while for $\sigma \neq 1/2$ it has first-order approximation.
	
	\begin{thm}\label{t-2}
		The difference scheme \eqref{4.9}--\eqref{4.12} is unconditionally stable for $\sigma \geq 1/2$. Under these conditions, the solution satisfies the a priori estimate:
		\begin{equation}\label{4.13}
			\|y^{n+1}\|^2_D + \sum_{k=1}^{m} a_k \|w_k^{n+1}\|^2 \leq \|u_0(x)\|^2_D,
			\quad n = 0,1,\ldots.
		\end{equation} 
	\end{thm}
	
	\begin{pf}
		Considering that
		\[
		y^{n+\sigma} = \left(\sigma - \frac{1}{2}\right) \tau \frac{y^{n+1} - y^{n}}{\tau} + \frac{y^{n+1} + y^{n}}{2},
		\]
		we rewrite \eqref{4.9}, \eqref{4.10} as
		\begin{equation}\label{4.14}
			\begin{split}
				\left(1 + r^h_m + \sum_{k=1}^{m} \frac{a_k}{\lambda_k}\right)
				\left(\frac{y^{n+1} - y^{n}}{\tau}, z\right) & +
				\left(\sigma - \frac{1}{2}\right) \tau \left(D \nabla \frac{y^{n+1} - y^{n}}{\tau}, \nabla z\right) \\
				& - 
				\sum_{k=1}^{m} \frac{a_k}{\lambda_k} \left(\frac{w_k^{n+1}-w_k^{n}}{\tau}, z\right) + \left(D \nabla \frac{y^{n+1} + y^{n}}{2}, \nabla z\right) = 0,
			\end{split}
		\end{equation} 
		\begin{equation}\label{4.15}
			\begin{split}    
				\left(1 + \left(\sigma - \frac{1}{2}\right) \lambda_k \tau \right) &
				\left(\frac{w_k^{n+1} - w_k^{n}}{\tau}, z_k\right) + 
				\lambda_k \left(\frac{w_k^{n+1} + w_k^{n}}{2}, z_k\right) \\
				& -
				\left(\frac{y^{n+1} - y^{n}}{\tau}, z_k\right) = 0,
				\quad k = 1,2,\ldots,m, 
				\quad n = 0,1,\ldots.
			\end{split}    
		\end{equation} 
		Setting
		\[
		z = 2(y^{n+1} - y^{n}),
		\quad z_k = 2\frac{a_k}{\lambda_k}(w_k^{n+1} - w_k^{n}),
		\quad k = 1,2,\ldots,m,
		\]
		from \eqref{4.14}, \eqref{4.15} with the constraint $\sigma \geq 1/2$, we obtain
		\[
		\begin{split}
			\left(D \nabla y^{n+1}, D \nabla y^{n+1}\right) & - \left(D \nabla y^{n}, D \nabla y^{n}\right) 
			+ \frac{2}{\tau} \sum_{k=1}^{m} \frac{a_k}{\lambda_k} \left(y^{n+1}-y^{n}, y^{n+1} - y^{n}\right) \\
			& - \frac{2}{\tau} \sum_{k=1}^{m} \frac{a_k}{\lambda_k} \left(w_k^{n+1}-w_k^{n}, y^{n+1} - y^{n}\right) \leq 0,
		\end{split}
		\]
		\[
		a_k(w_k^{n+1}, w_k^{n+1}) - a_k(w_k^{n}, w_k^{n}) -
		\frac{2}{\tau} \frac{a_k}{\lambda_k} \left(w_k^{n+1}-w_k^{n}, y^{n+1} - y^{n}\right) \leq 0,
		\quad k = 1,2,\ldots,m.
		\]
		Adding these inequalities, we obtain the solution estimate when advancing to the new time level:
		\[
		\left(D \nabla y^{n+1}, D \nabla y^{n+1}\right) + \sum_{k=1}^{m} a_k(w_k^{n+1}, w_k^{n+1}) \leq \left(D \nabla y^{n}, D \nabla y^{n}\right) + \sum_{k=1}^{m} a_k(w_k^{n}, w_k^{n}),
		\]
		from which \eqref{4.13} follows.
		\qed
	\end{pf}

	A brief comment is appropriate regarding the computational implementation of the considered two-level scheme \eqref{4.9}--\eqref{4.12}. In our case, the same finite element basis is used for both $y(x,t)$ and $w_k(x,t)$, $k = 1,2,\ldots,m$. This allows us to express $w_k^{n+1}$ from \eqref{4.10} as follows:
	\begin{equation}\label{4.16}
		(w_k^{n+1}, z) = \frac{1}{1 + \sigma \lambda_k \tau} (y^{n+1}, z) + (f_k^{n}, z),
	\end{equation}  
	where 
	\[
	(f_k^{n}, z) = \frac{1}{1 + \sigma \lambda_k \tau} \Big(\big(1 - (1-\sigma) \lambda_k \tau\big) (w_k^{n}, z) - (y^{n},z)\Big),
	\quad k = 1,2,\ldots,m.
	\]
	Substituting into \eqref{4.9}, we obtain
	\begin{equation}\label{4.17}
		(1 + r^h_m + \sigma \tau \alpha) (y^{n+1}, z) + \sigma \tau \left(D \nabla y^{n+1}, \nabla z\right) = (f^n, z),
	\end{equation} 
	for finding $y^{n+1}$. Here we use the notations
	\[
	\alpha = \sum_{k=1}^{m} \frac{a_k}{1 + \sigma \lambda_k \tau},
	\quad f^n = (1 + r^h_m) (y^{n},z) - (1-\sigma)\tau \left(D \nabla y^{n}, \nabla z\right) + \sum_{k=1}^{m} \frac{a_k}{\lambda_k} (y^{n} - w_k^{n} + f_k^n, z).
	\] 
	The transition to the new time level $n+1$ is achieved by solving the standard problem (\ref{4.17}) for $y^{n+1}$ and computing the auxiliary quantities $w_k^{n+1}$, $k = 1,2,\ldots,m$, according to (\ref{4.16}).
	The computational complexity of the approximate solution to the considered nonlocal problem is not much greater than that of a standard local problem for a parabolic equation. It only requires additionally solving $m$ simple auxiliary local evolutionary problems using explicit procedures.
	
	\section{Numerical experiments }\label{sec:5} 
	
	The capabilities of the described computational algorithm were evaluated through numerical solution of test problems. We present results for computing the effective diffusion coefficient, approximating the memory kernel, and applying these results to the homogenized nonstationary diffusion problem in a periodic medium.
	
	\begin{figure}[ht]
		\centering
		\includegraphics[width=1\linewidth]{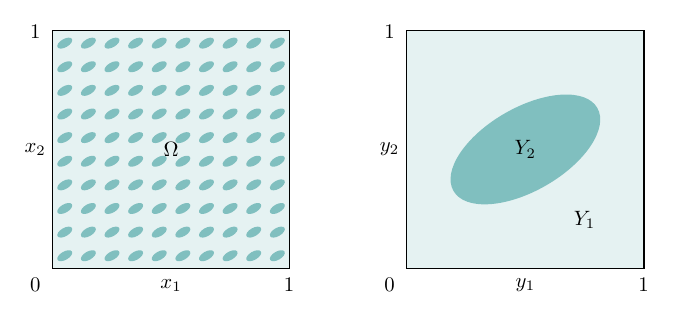}
		\caption{Computational domain $\Omega$ (left) and periodicity cell (right) for the test problem.}
		\label{f-2}
	\end{figure}
	
	\subsection{Test problem}
	
	We illustrate the main elements of the computational technology by solving problem \eqref{2.1}--\eqref{2.3} approximately, where the computational domain $\Omega$ is a unit square and the inclusion cross-section is an ellipse with 2:1 aspect ratio rotated by $30^\circ$, occupying about one quarter ($0.08\pi$) of the cell (Fig.~\ref{f-2}). 
	For the diffusion coefficient \eqref{2.4}, we set
	\[
	d(y) =
	\begin{cases}
		1, & y \in Y_1, \\
		\varepsilon^2, & y \in Y_2,
	\end{cases} 
	\]
	i.e., $d_1 = d_2 = 1$. The initial condition \eqref{2.2} is specified as
	\[
	u_0(x) = \frac{4}{1 + \exp(-100(x_1-0.5))} x_1(1-x_1) \sin(\pi x_2).
	\]
	
	\begin{figure}[ht]
		\centering
		\includegraphics[width=0.48\linewidth]{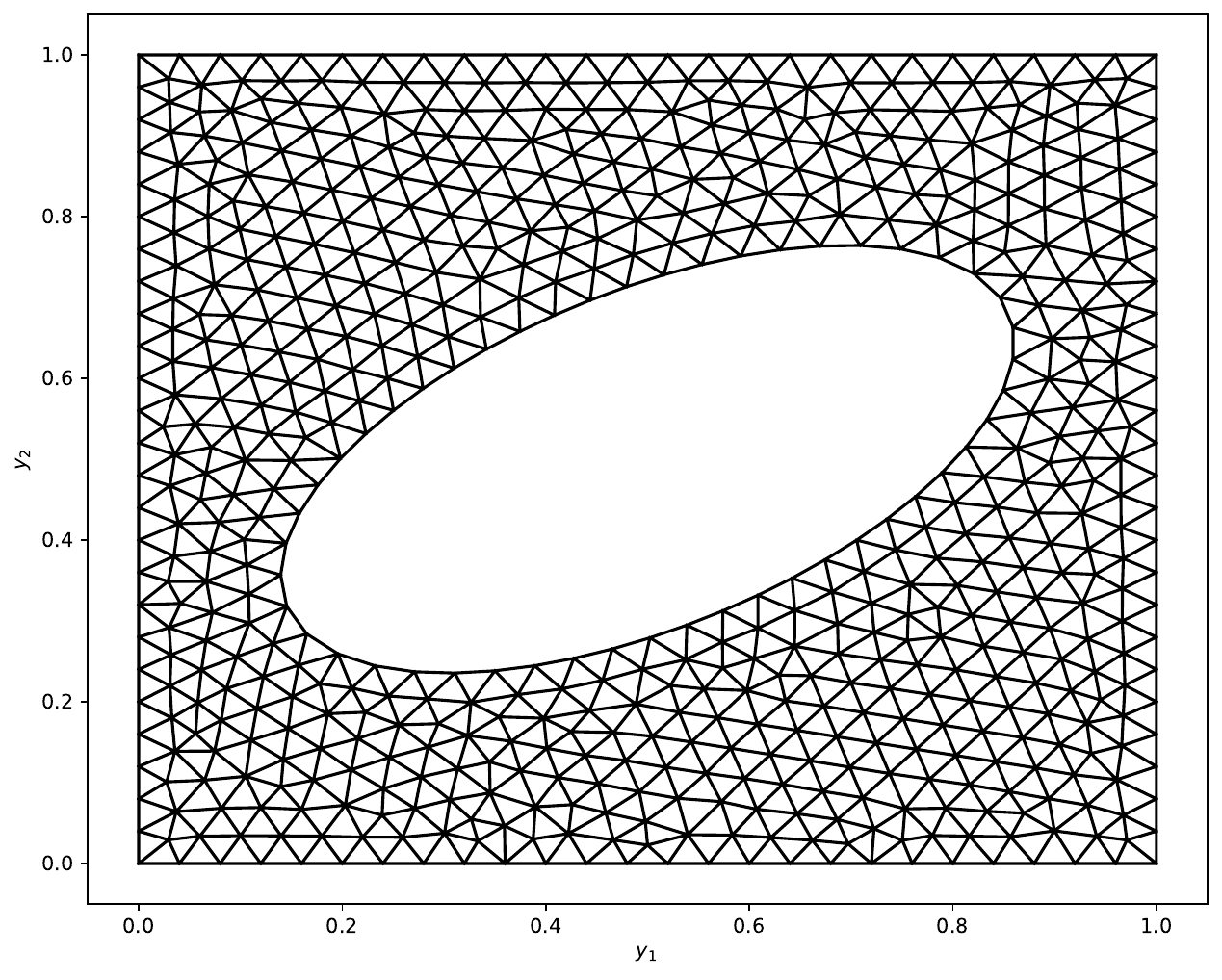}
		\includegraphics[width=0.48\linewidth]{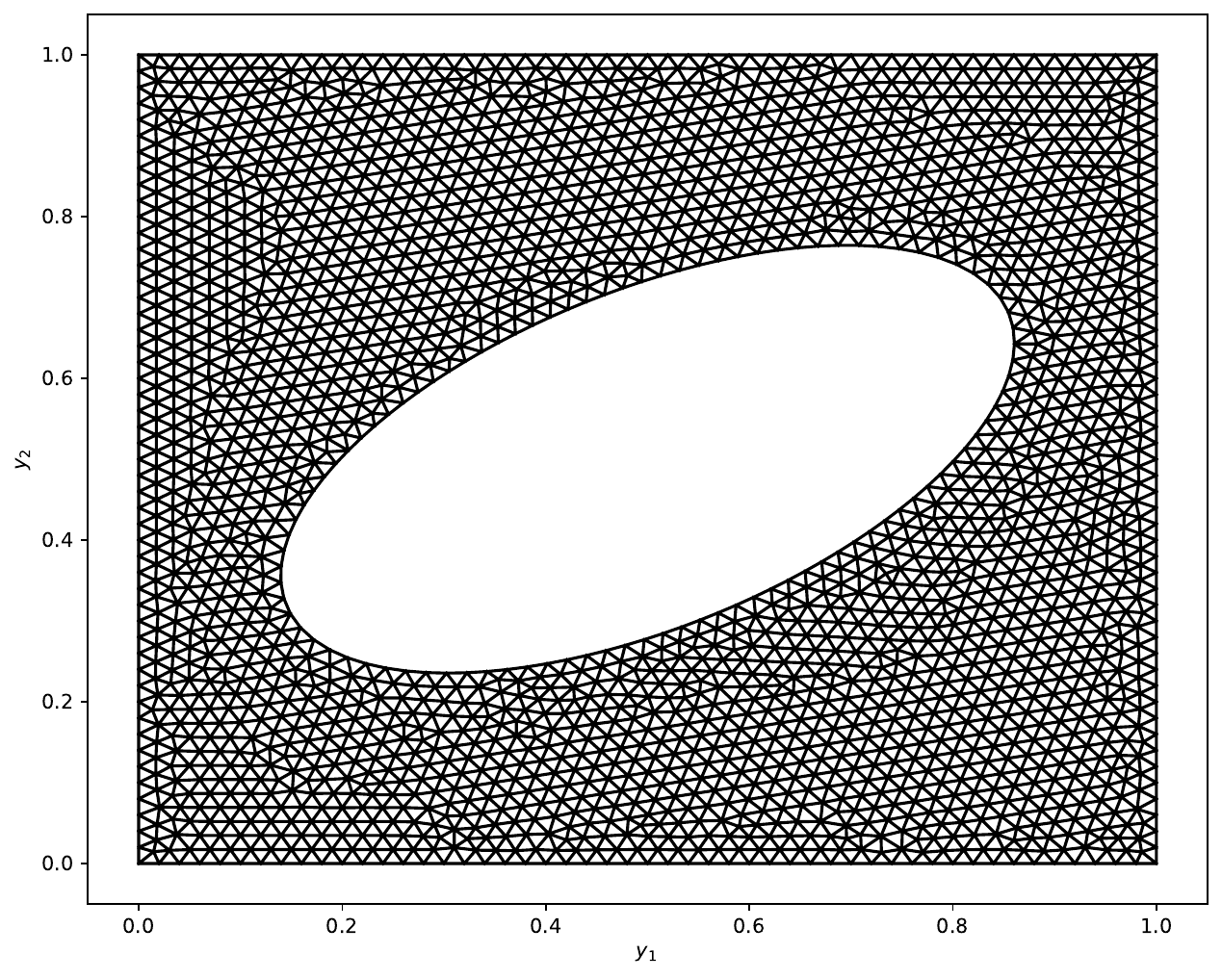}
		\caption{Grid 1 (left) in $Y_1$ and Grid 2 (right).}
		\label{f-5.1}
	\end{figure}
	
	\begin{figure}[ht]
		\centering
		\includegraphics[width=0.48\linewidth]{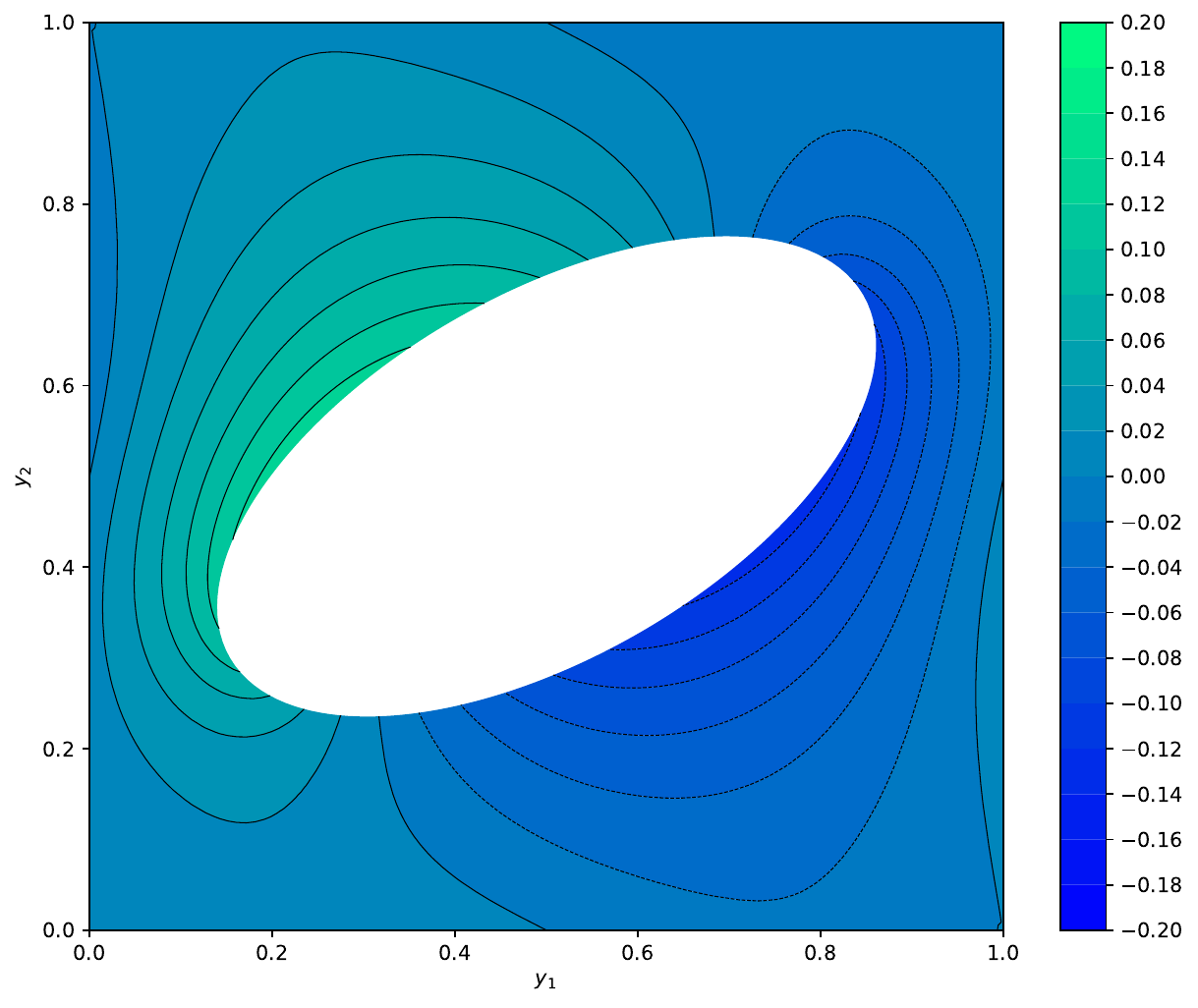}
		\includegraphics[width=0.48\linewidth]{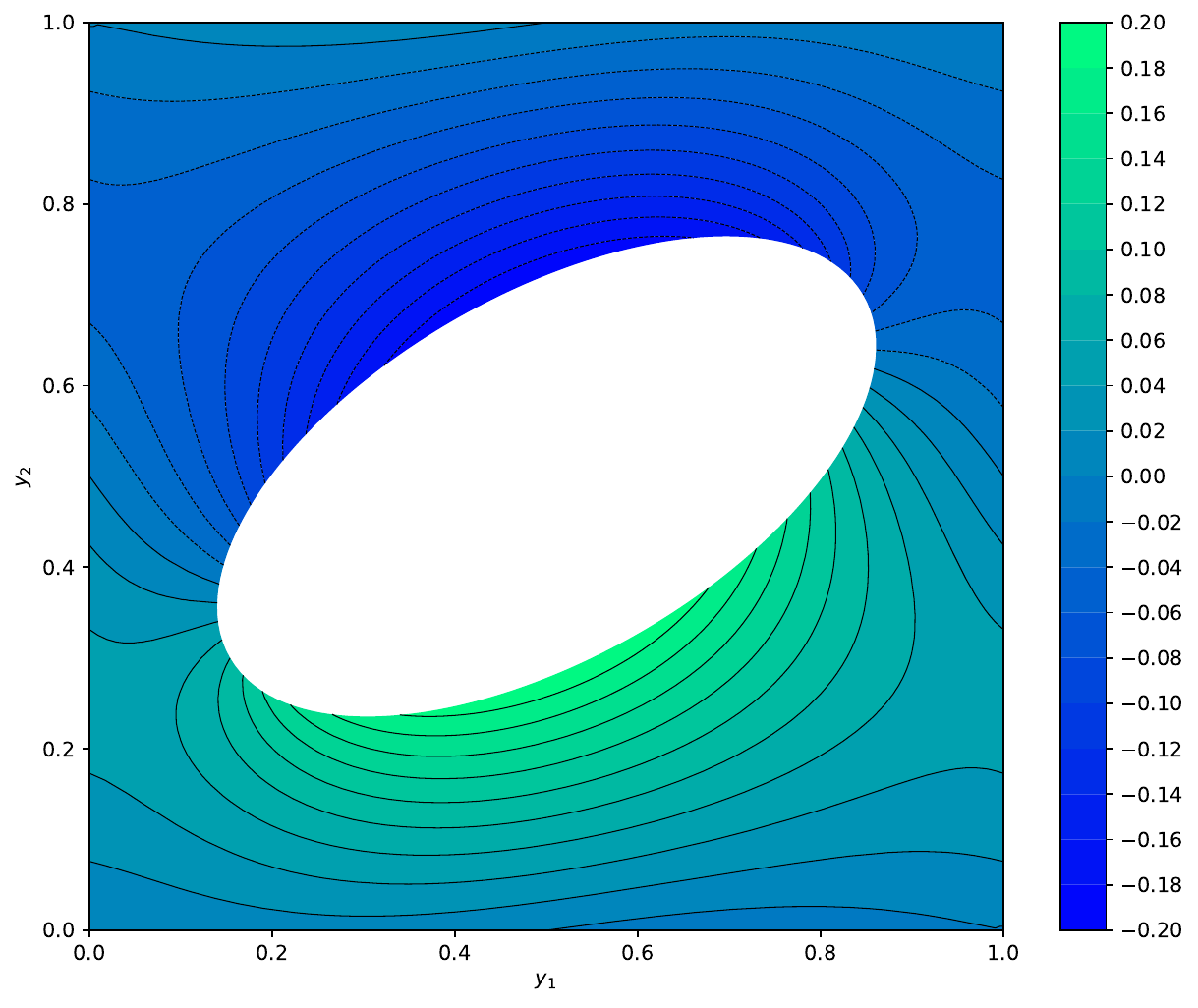}
		\caption{Auxiliary functions: $\theta_1^h(y)$ (left) and $\theta_2^h(y)$ (right).}
		\label{f-5.2}
	\end{figure}
	
	\subsection{Diffusion tensor}
	
	To find the effective diffusion tensor $D$ in equation \eqref{2.5}, we solve the boundary-value problem on the periodicity cell $Y_1$. We present calculation results on three grids: Grid 1 (648 nodes), Grid 2 (2375 nodes), and Grid 3 (9075 nodes), generated using Gmsh \cite{geuzaine2009gmsh} software. Grids 1 and 2 are shown in Fig.~\ref{f-5.1}.
	First, according to \eqref{4.1}, we find the auxiliary functions $\theta_i^h(y)$, $i = 1,2$. These are shown in Fig.~\ref{f-5.2} using Grid 3 with quadratic Lagrange finite elements. The software implementation utilizes the FEniCS computing platform \cite{alnaes2015fenics}.
	
	The diffusion tensor is calculated from the auxiliary functions $\theta_i^h(y)$, $i = 1,2$ according to \eqref{4.2}. Numerical results using different grids and finite element analysis (FEA) with linear ($p=1$) and quadratic ($p=2$) elements are presented in Table~1. The symmetry property of the diffusion tensor is maintained at the discrete level. We note that even on coarse grids, the accuracy of individual tensor components is high.
	
	\begin{table}
		\label{tab-1}
		\centering
		\caption{Effective diffusion tensor}
		\vspace{1ex}
		\begin{tabular}{lcc}
			\toprule
			FEA & $p=1$ & $p=2$ \\
			\midrule
			Grid 1 & 
			$\begin{pNiceMatrix}[columns-width=auto,r]
				0.86027695 & -0.08233093 \\
				-0.08233093 & 0.97912460
			\end{pNiceMatrix}$ &
			$\begin{pNiceMatrix}[columns-width=auto,r]
				0.86027824 & -0.08198406 \\
				-0.08198406 & 0.98004661
			\end{pNiceMatrix}$ \\ [2ex]
			
			Grid 2 &
			$\begin{pNiceMatrix}[columns-width=auto,r]
				0.85961946 & -0.08227680 \\
				-0.08227680 & 0.97969683
			\end{pNiceMatrix}$ &
			$\begin{pNiceMatrix}[columns-width=auto,r]
				0.85962460 & -0.08218294 \\
				-0.08218294 & 0.97993062 \\
			\end{pNiceMatrix}$ \\ [2ex]
			
			Grid 3 &
			$\begin{pNiceMatrix}[columns-width=auto,r]
				0.85945517 & -0.08225859 \\
				-0.08225859 & 0.97984380
			\end{pNiceMatrix}$ &
			$\begin{pNiceMatrix}[columns-width=auto,r]
				0.85945621 & -0.08223535 \\
				-0.08223535 & 0.97990134
			\end{pNiceMatrix}$ \\
			\bottomrule
		\end{tabular}
	\end{table}
	
	\begin{figure}[ht]
		\centering
		\includegraphics[width=0.48\linewidth]{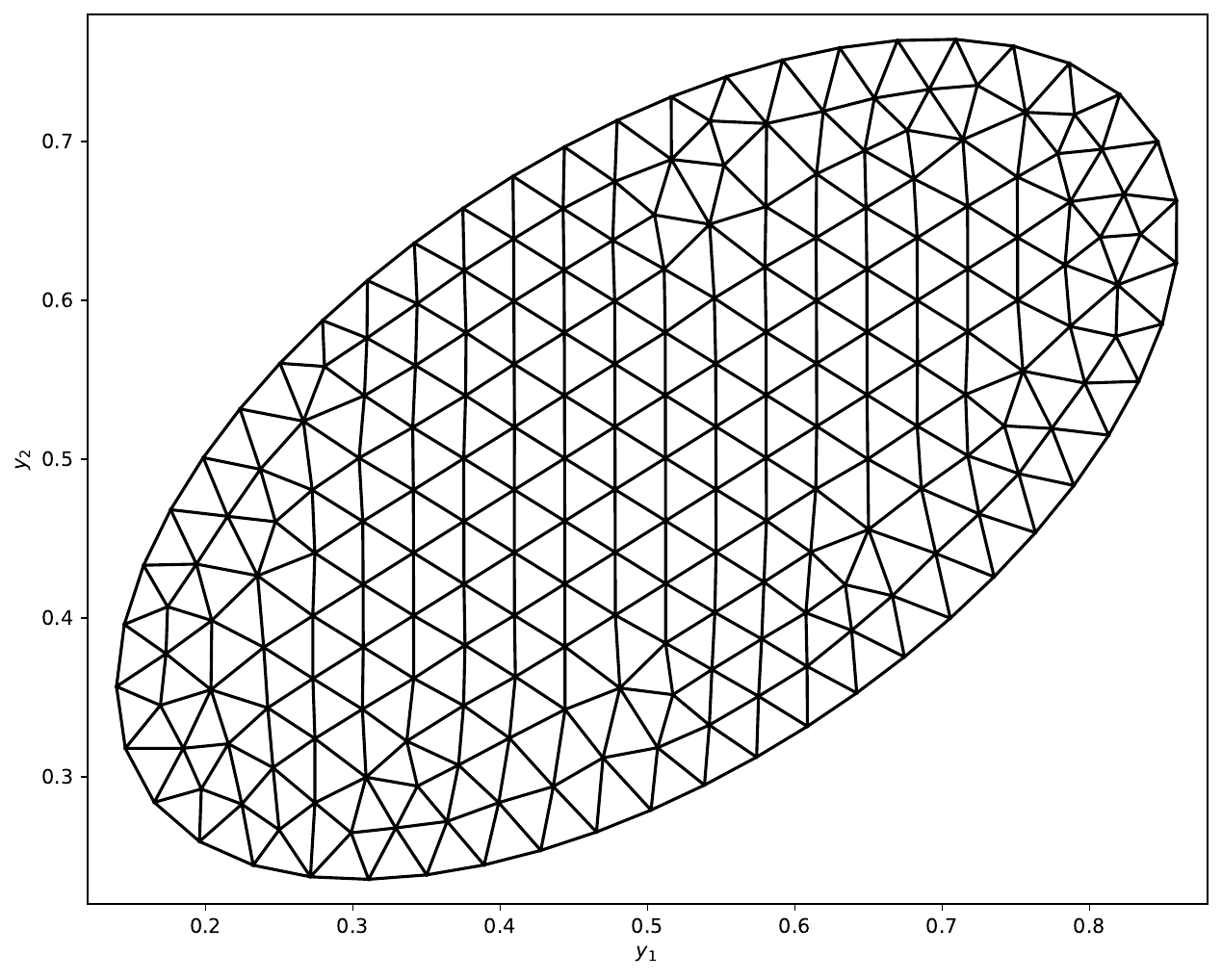}
		\includegraphics[width=0.48\linewidth]{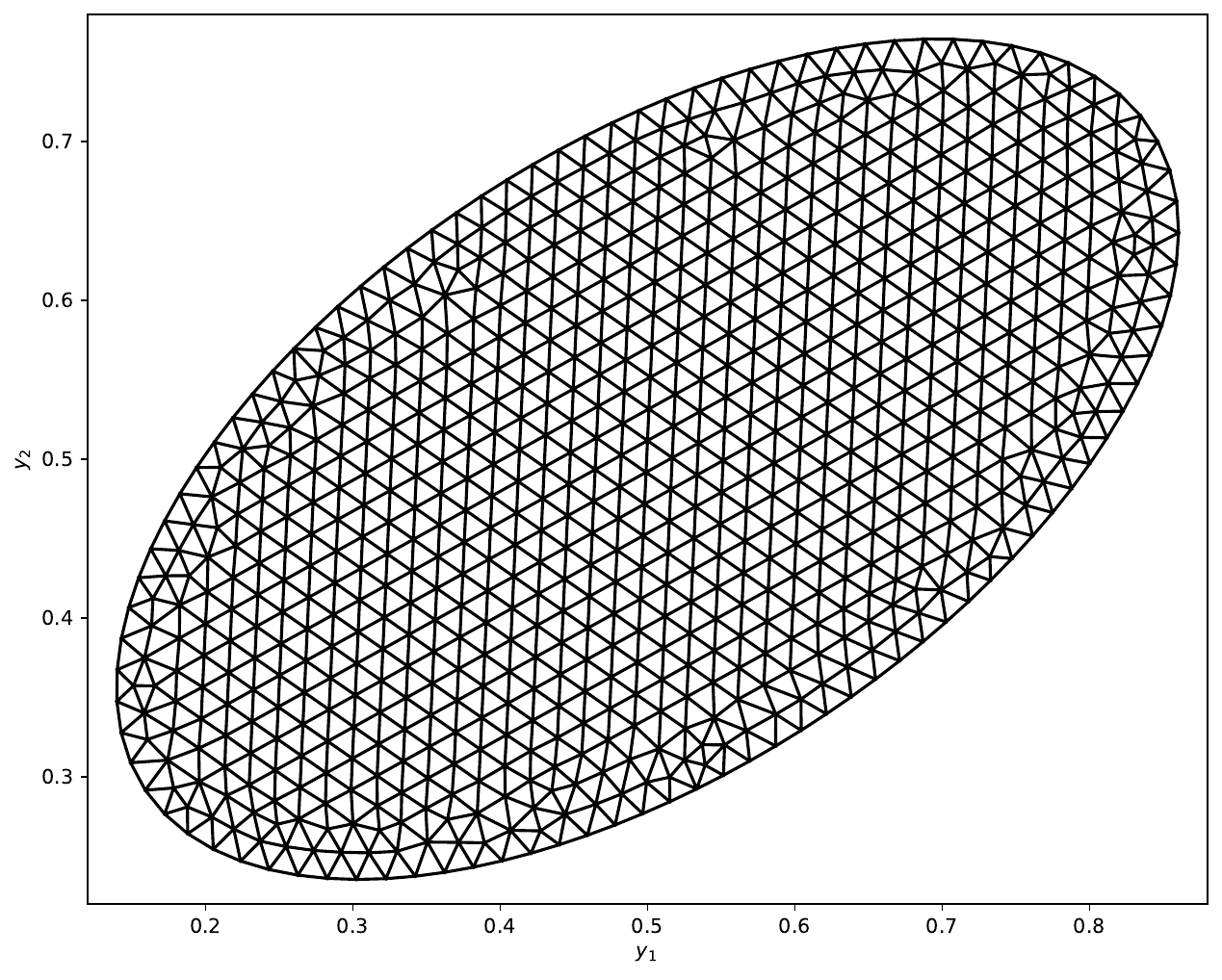}
		\caption{Grid 1 (left) in $Y_2$  and Grid 2 (right).}
		\label{f-5.3}
	\end{figure}
	
	\subsection{Memory kernel}
	
	\begin{figure}[ht]
		\centering
		\includegraphics[width=0.48\linewidth]{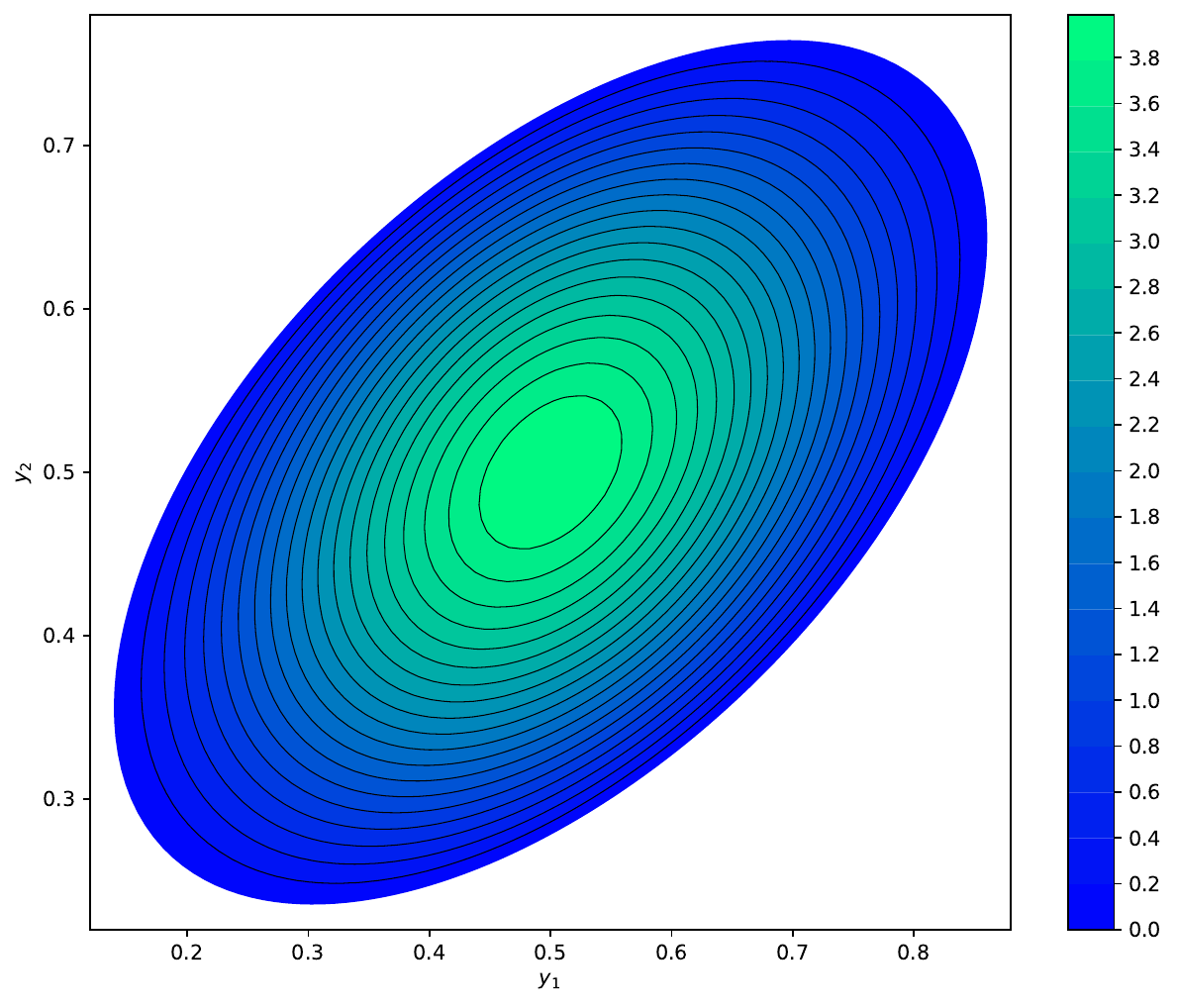}
		\includegraphics[width=0.48\linewidth]{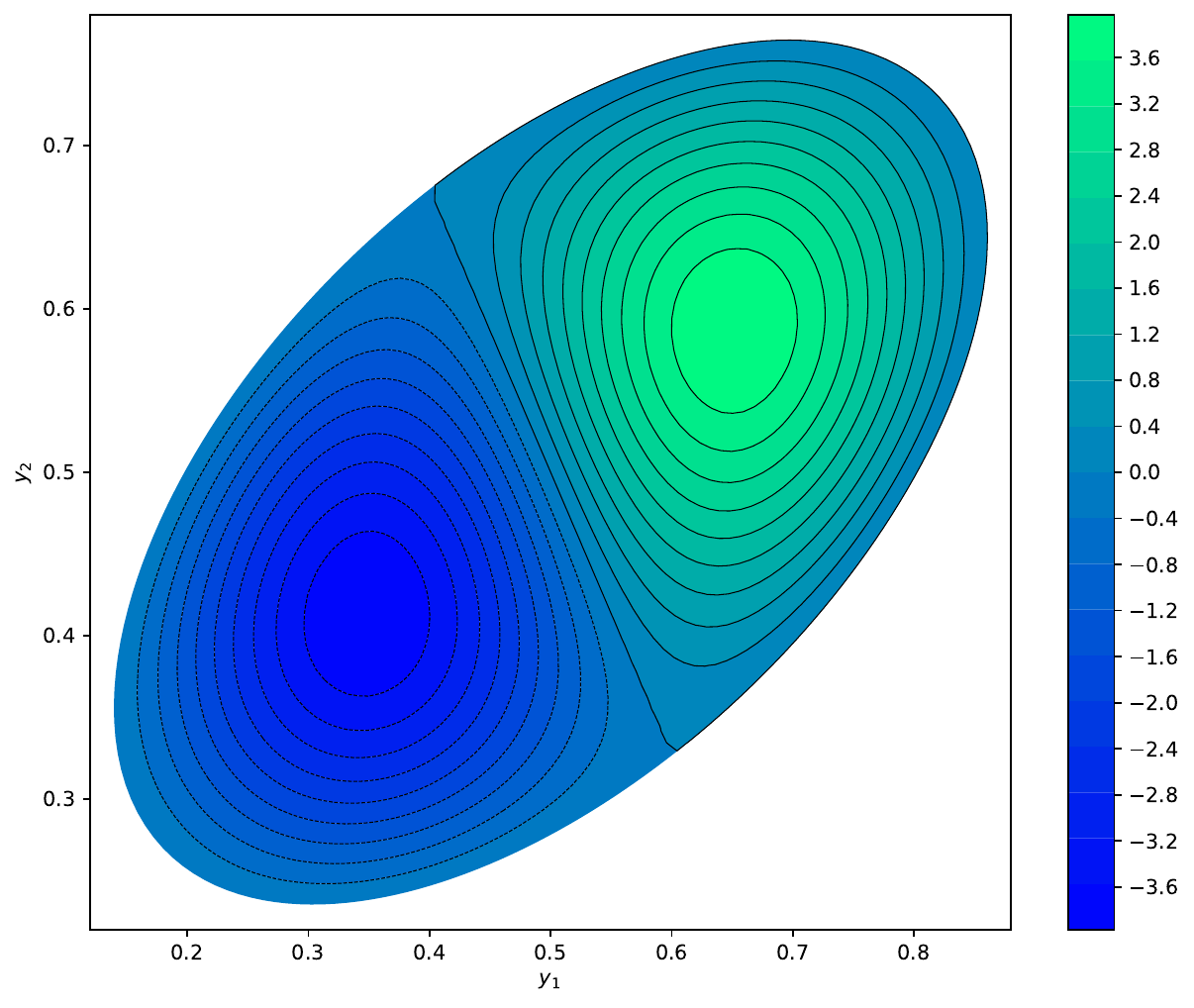} \\
		\includegraphics[width=0.48\linewidth]{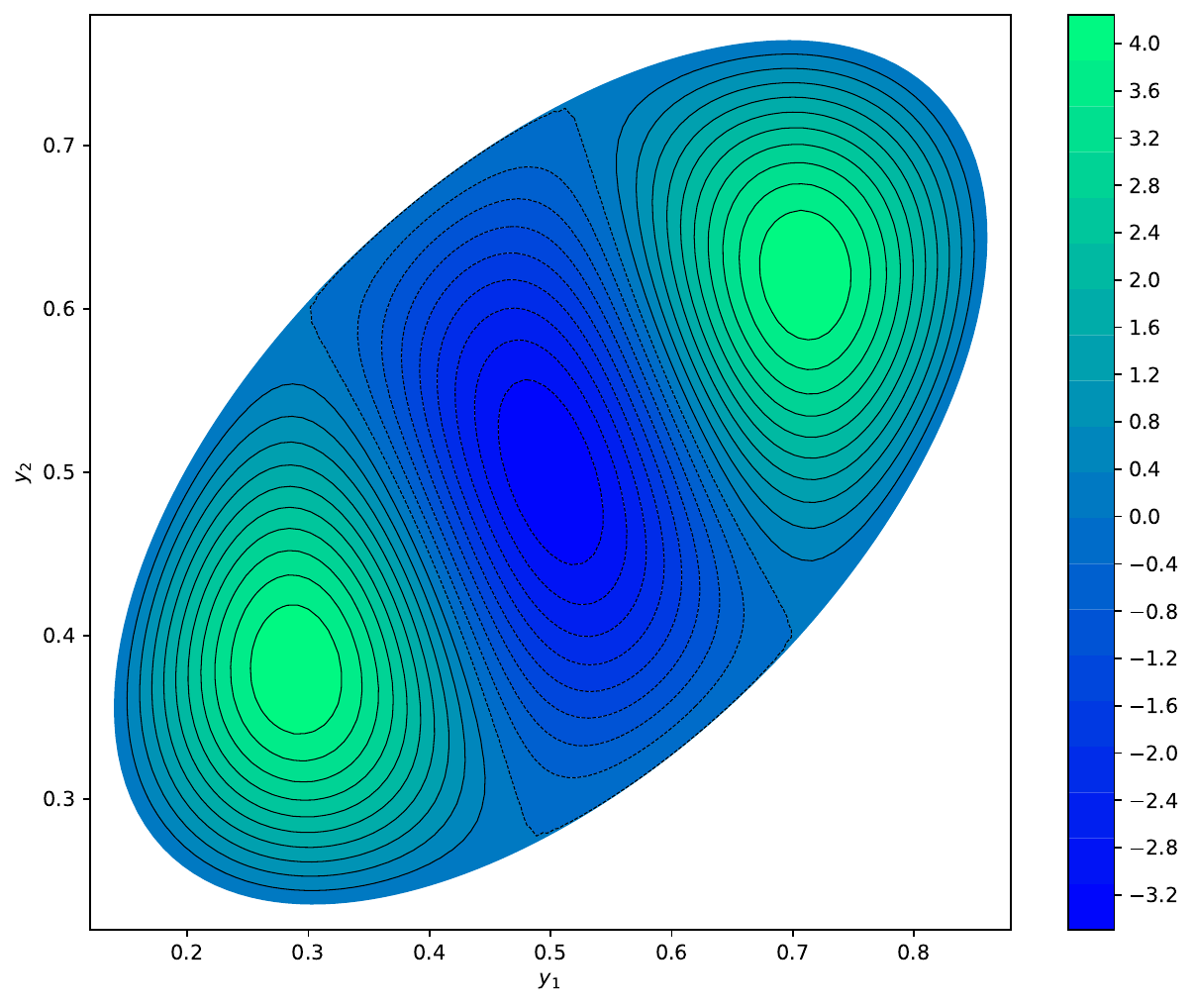}
		\includegraphics[width=0.48\linewidth]{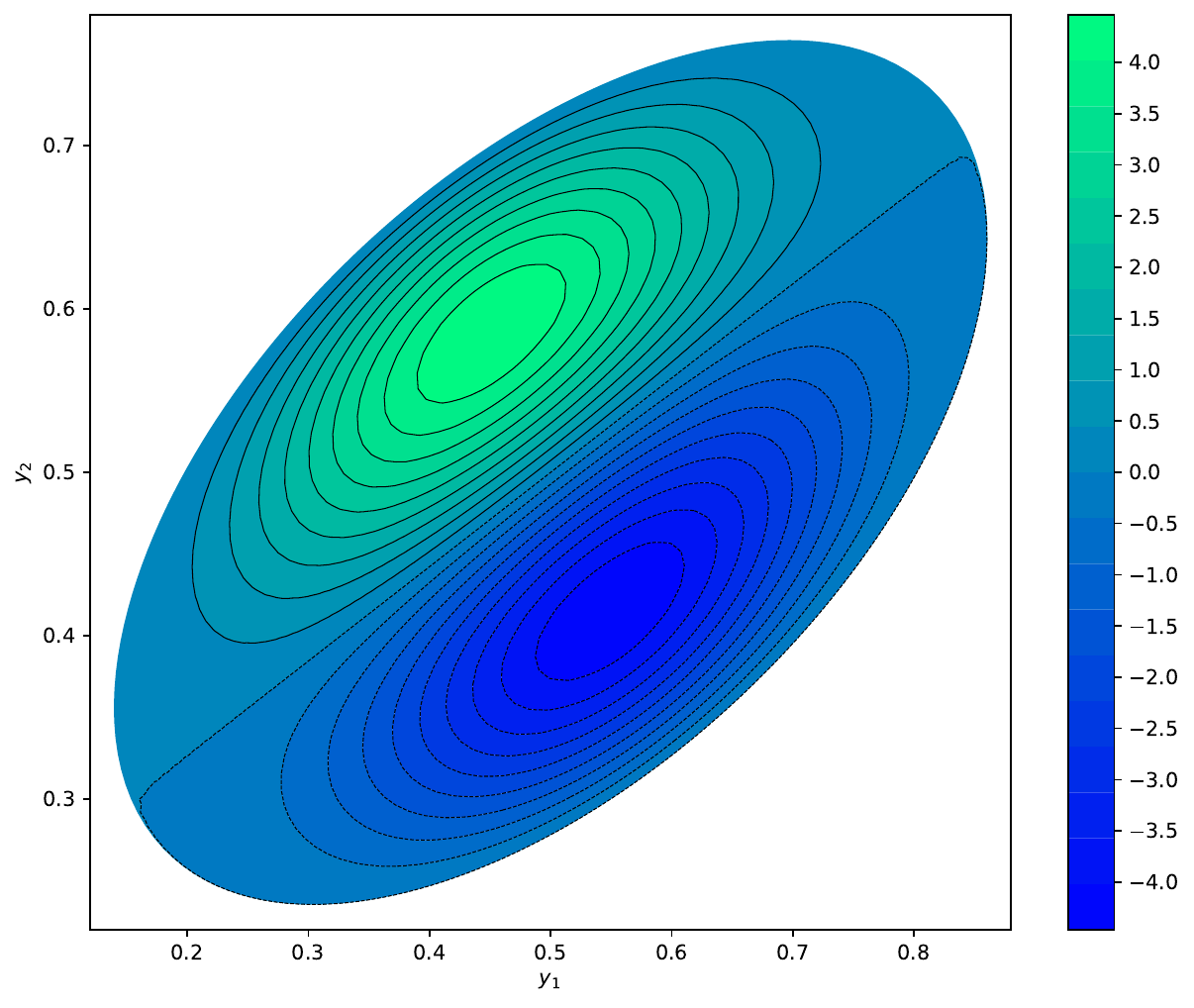} \\
		\caption{Eigenfunctions: top row --- $\varphi^h_1(y)$ (left) and $\varphi^h_2(y)$ (right), bottom row --- $\varphi^h_3(y)$ (left) and $\varphi^h_4(y)$ (right).}
		\label{f-5.4}
	\end{figure}
	
	First, we solve the spectral problem \eqref{4.3} on the part $Y_2$ of the periodicity cell. We use computational grids: Grid 1 (224 nodes), Grid 2 (797 nodes), and Grid 3 (3055 nodes). Grids 1 and 2 are shown in Fig.~\ref{f-5.3}.
	
	The numerical results of solving the spectral problem \eqref{4.3} on different grids using first- and second-order finite elements are presented in Table~2. We observe convergence of eigenvalues as the grid is refined and the polynomial order $p$ is increased. The first four eigenfunctions obtained using the finest grid and quadratic finite elements are shown in Fig.~\ref{f-5.4}.
	
	\begin{table}[htp]
		\label{tab-2}
		\caption{First 10 eigenvalues}
		\vspace{1ex}
		\centering
		\begin{tabular}{l c ll c ll c ll}
			\toprule
			& \multicolumn{2}{c}{Grid 1} && \multicolumn{2}{c}{Grid 2} && \multicolumn{2}{c}{Grid 3}\\
			\cline{2-3} \cline{5-6} \cline{8-9} 
			$k$ & $p=1$ & $p=2$ && $p=1$ & $p=2$ && $p=1$ & $p=2$ \\
			\midrule
			1 & 8.9942e+01 & 8.9369e+01 && 8.9369e+01 & 8.9221e+01 && 8.9221e+01 & 8.9184e+01 \\
			2 & 1.5925e+02 & 1.5737e+02 && 1.5751e+02 & 1.5701e+02 && 1.5704e+02 & 1.5692e+02 \\
			3 & 2.5686e+02 & 2.5169e+02 && 2.5228e+02 & 2.5095e+02 && 2.5111e+02 & 2.5078e+02 \\
			4 & 3.0198e+02 & 2.9401e+02 && 2.9558e+02 & 2.9357e+02 && 2.9396e+02 & 2.9346e+02 \\
			5 & 3.8516e+02 & 3.7365e+02 && 3.7539e+02 & 3.7235e+02 && 3.7281e+02 & 3.7205e+02 \\
			6 & 4.1396e+02 & 3.9907e+02 && 4.0209e+02 & 3.9833e+02 && 3.9910e+02 & 3.9817e+02 \\
			7 & 5.4730e+02 & 5.2394e+02 && 5.2794e+02 & 5.2183e+02 && 5.2288e+02 & 5.2134e+02 \\
			8 & 5.5296e+02 & 5.2714e+02 && 5.3249e+02 & 5.2594e+02 && 5.2731e+02 & 5.2568e+02 \\
			9 & 6.6068e+02 & 6.2370e+02 && 6.3188e+02 & 6.2247e+02 && 6.2458e+02 & 6.2224e+02 \\
			10 & 7.2158e+02 & 6.7946e+02 && 6.8853e+02 & 6.7755e+02 && 6.7989e+02 & 6.7716e+02 \\
			\bottomrule
		\end{tabular}
	\end{table}
	
	\clearpage
	
	\begin{figure}[ht]
		\centering
		\includegraphics[width=0.48\linewidth]{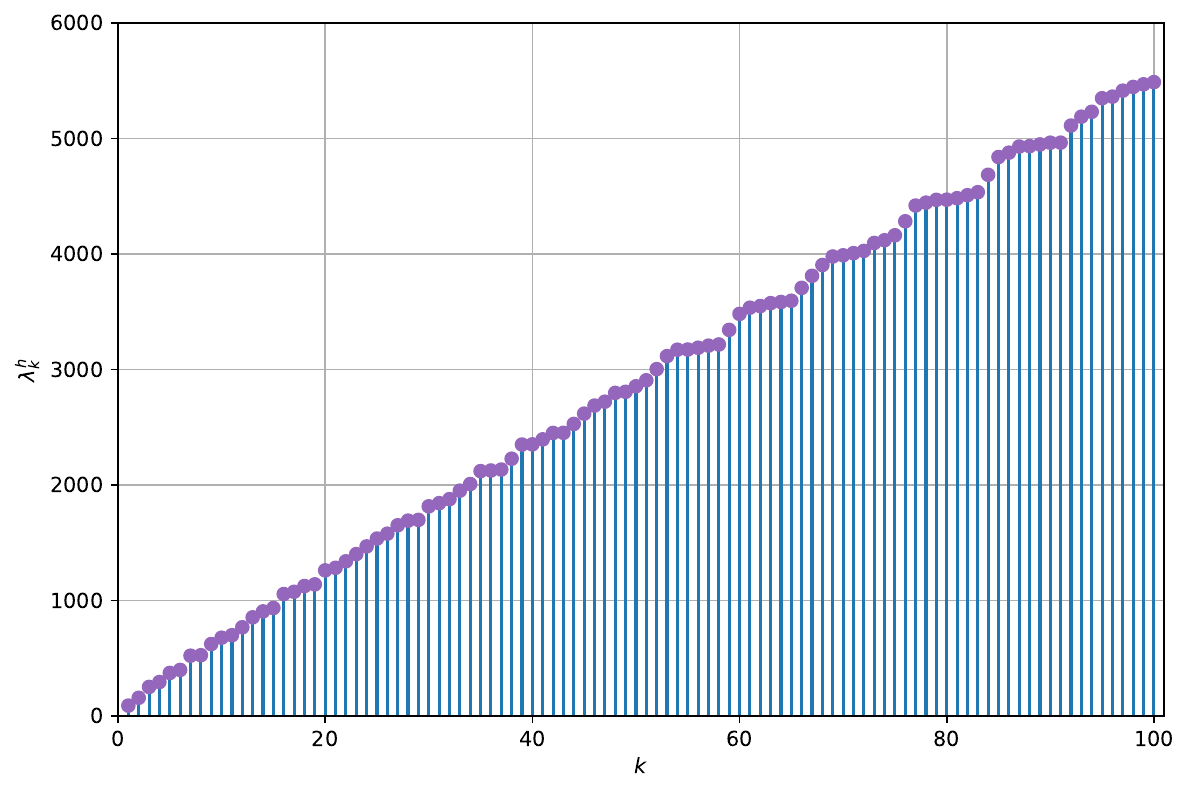}
		\includegraphics[width=0.48\linewidth]{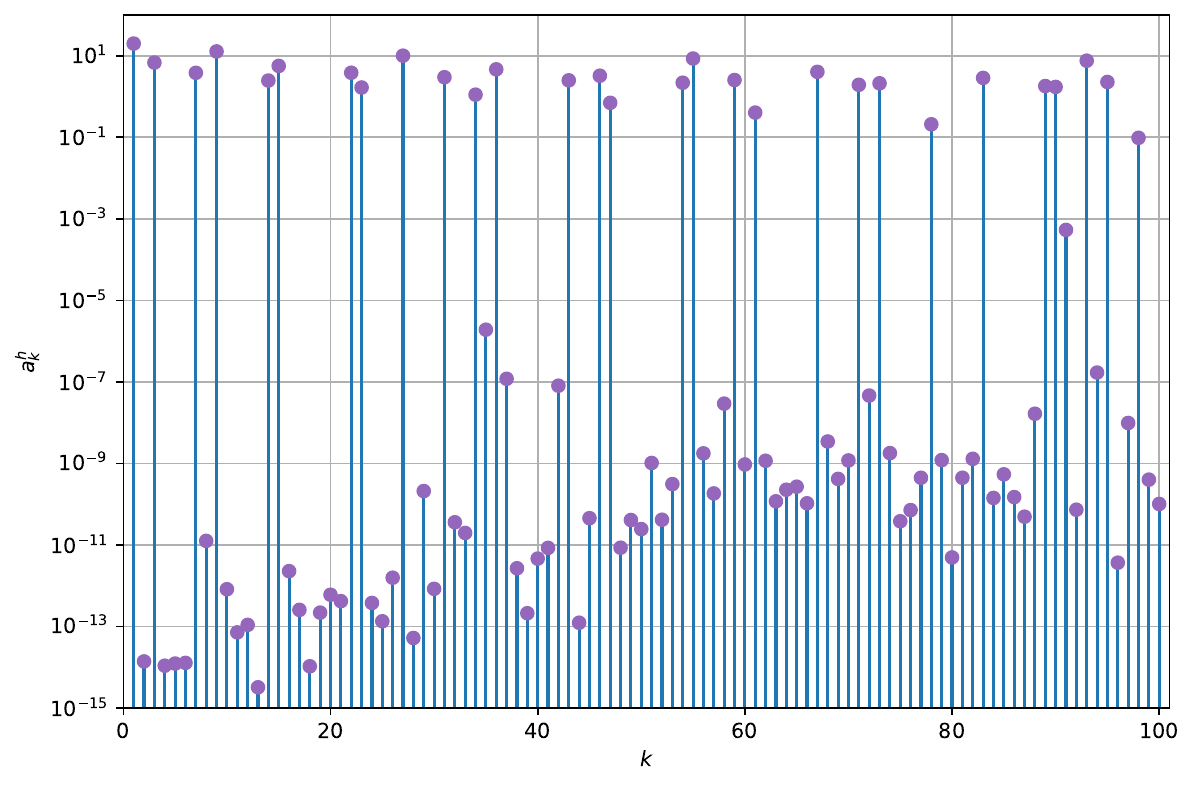}
		\caption{Coefficients for the exponential sum approximation.}
		\label{f-5.5}
	\end{figure}
	
	\begin{table}
		\label{t-3}
		\centering
		\caption{Filtering of small coefficients}
		\vspace{1ex}
		\begin{tabular}{llllll}
			\toprule
			$\epsilon$ & $10^{-5}$ & $10^{-6}$ & $10^{-7}$ & $10^{-8}$ & $10^{-9}$ \\
			\midrule
			$m(\epsilon)$ & $30$ & $31$ & $33$ & $37$ & $46$ \\
			$\varrho_m(\epsilon)$ & 2.3951e-06 & 4.8900e-07 & 2.0085e-07 & 2.8173e-08 & 5.4933e-09 \\
			\bottomrule
		\end{tabular}
	\end{table}
	
\begin{figure}[ht]
		\centering
		\includegraphics[width=0.48\linewidth]{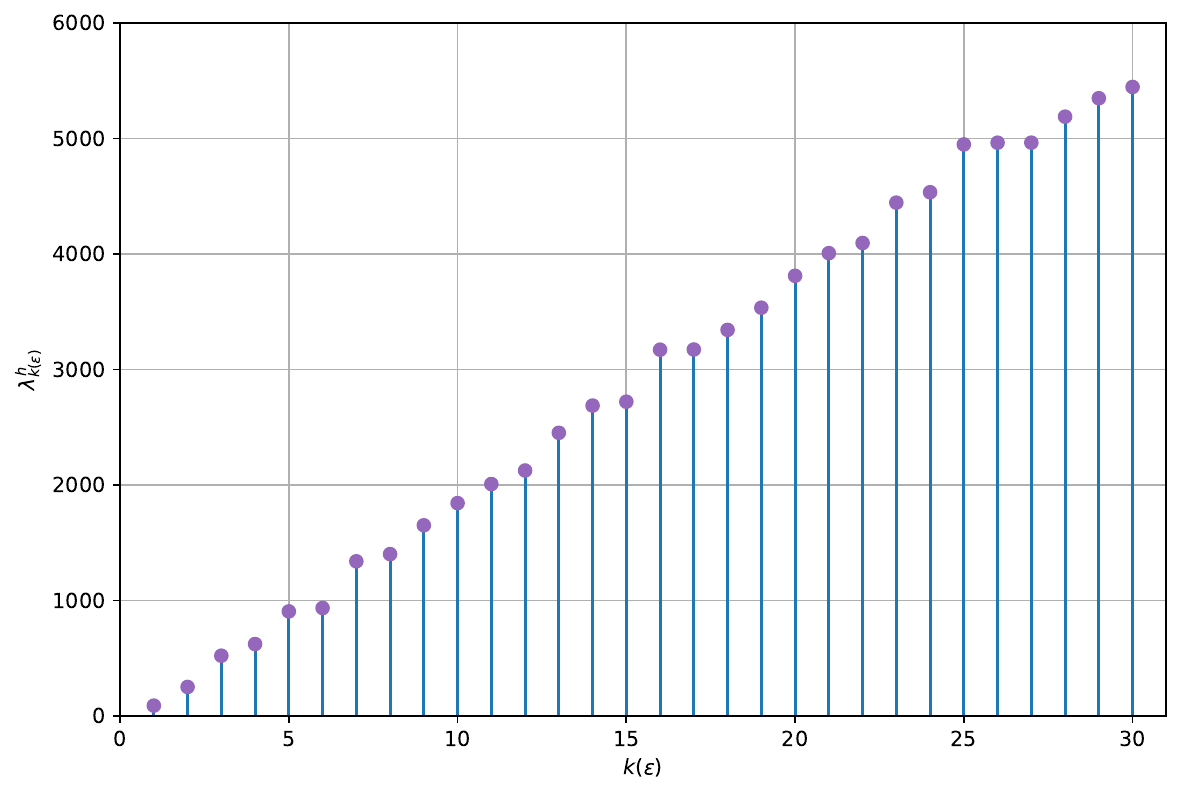}
		\includegraphics[width=0.48\linewidth]{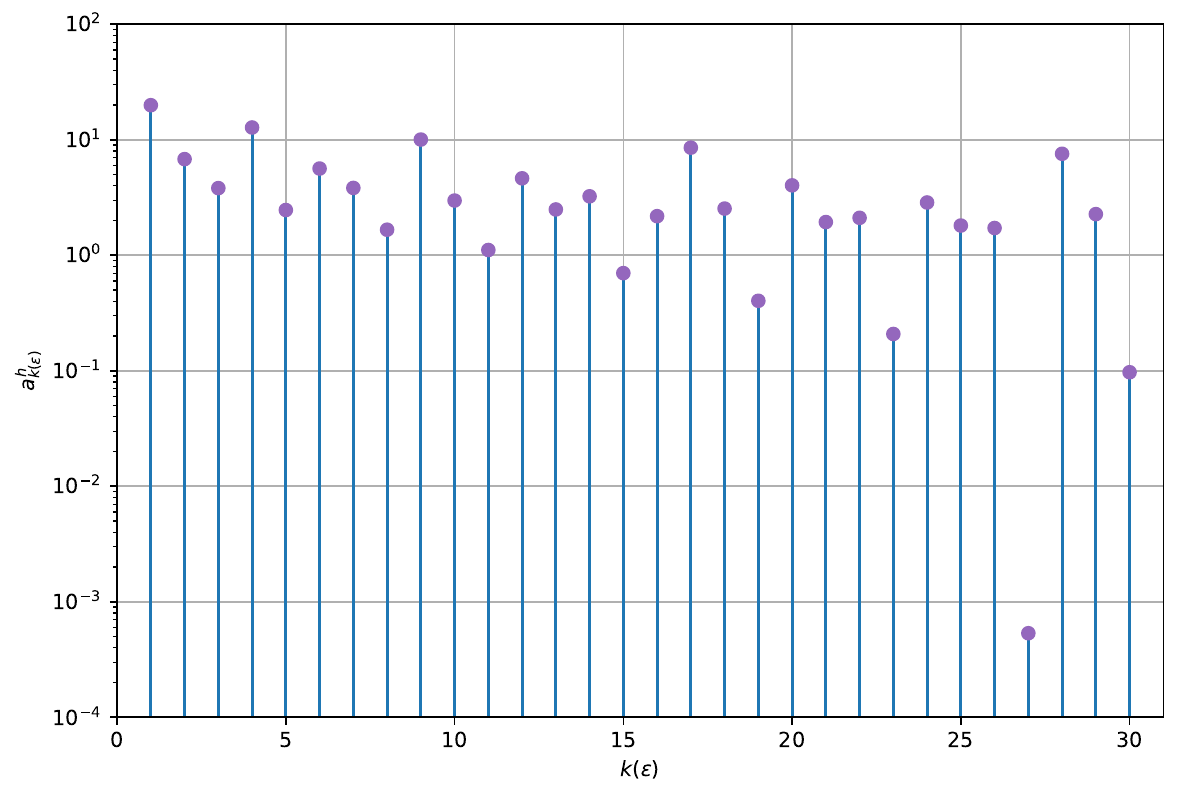}
		\caption{Filtered coefficients for the exponential sum approximation.}
		\label{f-5.6}
	\end{figure}

	For kernel approximation, we use the solution of the spectral problem on Grid 3 with $p=2$. The baseline data includes 100 eigenvalues ($m=100$ in approximation \eqref{4.3}). The eigenvalues $\lambda_k^h$ and coefficients $a_k^h$, $k=1,2,\dots,m$ are shown in Fig.~\ref{f-5.5}.
	We observe that only a relatively small number of eigenvalues have significant influence. There are many small coefficients $a_k^h$, $k=1,2,\dots,m$, which allows us to retain in representation \eqref{4.4} only those terms exceeding a threshold value $\epsilon$. 
	For $a_k < \epsilon$, the corresponding term in \eqref{4.4} is excluded, and after renumbering we have
	\[
	\chi_{m(\epsilon)}^h(t) = \sum_{k(\epsilon)=1}^{m(\epsilon)} a^h_{k(\epsilon)} \exp\left(-\lambda^h_{k(\epsilon)} t\right).
	\]
	The effect of this filtering can be estimated by $\varrho_m(\epsilon) = \chi_m^h(0) - \chi_{m(\epsilon)}^h(0)$. In our case, the filtering error estimates are presented in Table~3, where $m=100$ and $\chi_m^h(0) \approx 120.4433$.
	For this example, we can use about one-third of the eigenvalues without significant loss of accuracy. After filtering with $\epsilon=10^{-5}$, the eigenvalues $\lambda_k^h$ and corresponding coefficients $a_k^h$, $k=1,2,\dots,m$ are shown in Fig.~\ref{f-5.6}.
	
	\begin{figure}[ht]
		\centering
		\includegraphics[width=0.5\linewidth]{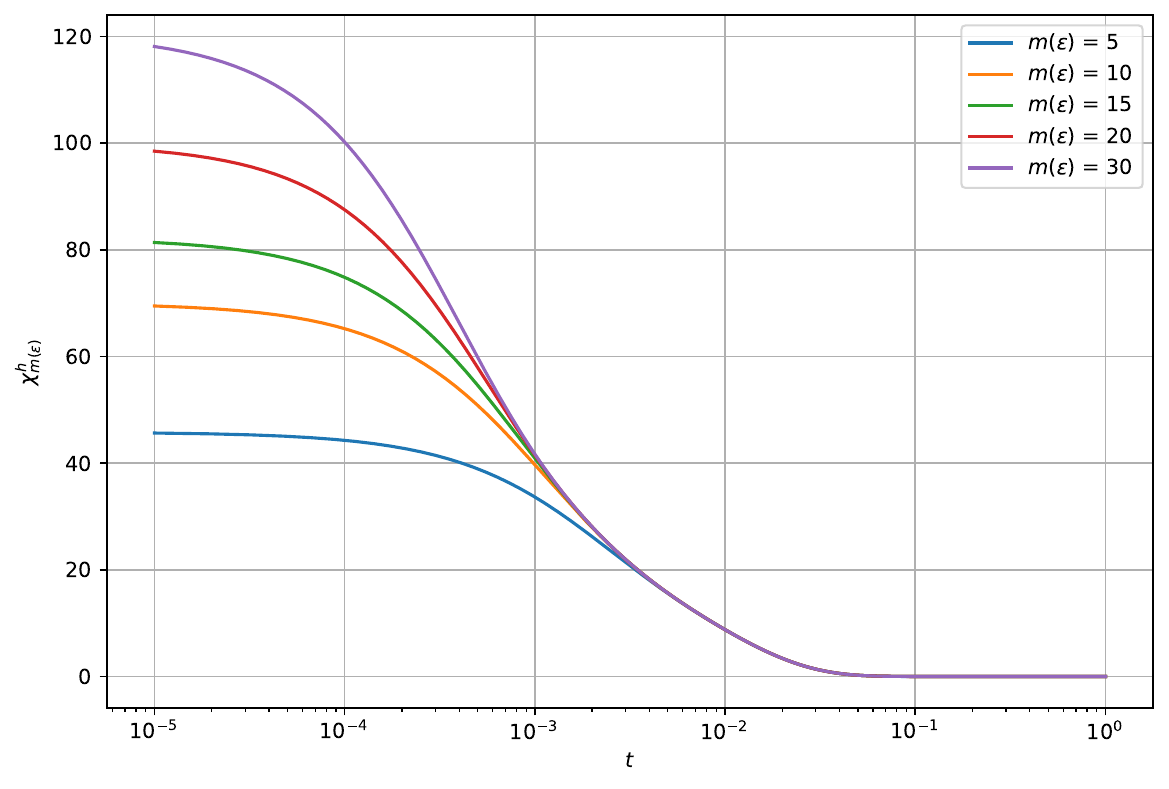}
		\caption{Approximation of the memory kernel by a sum of exponentials.}
		\label{f-5.7}
	\end{figure}
	
	The results of kernel approximation by exponential sums with filtering parameter $\epsilon=10^{-5}$ are presented in Fig.~\ref{f-5.7}. Only $m(\epsilon)$ leading eigenvalues exceeding $\epsilon=10^{-5}$ are used. The kernel has a singularity at $t\rightarrow 0$, so significant errors are expected for small $t$.
	
	\begin{table}
		\label{t-4}
		\centering
		\caption{Higher harmonics parameter}
		\vspace{1ex}
		\begin{tabular}{llllll}
			\toprule
			$m(\epsilon)$ & $5$ & $10$ & $15$ & $20$ & $30$ \\
			\midrule
			$m$ & $14$ & $31$ & $47$ & $67$ & $98$ \\
			$r_m^h$ & 0.054385 & 0.036592 & 0.031373 & 0.026059 & 0.021766 \\
			\bottomrule
		\end{tabular}
	\end{table}
	
\begin{figure}[ht]
		\centering
		\includegraphics[width=0.48\linewidth]{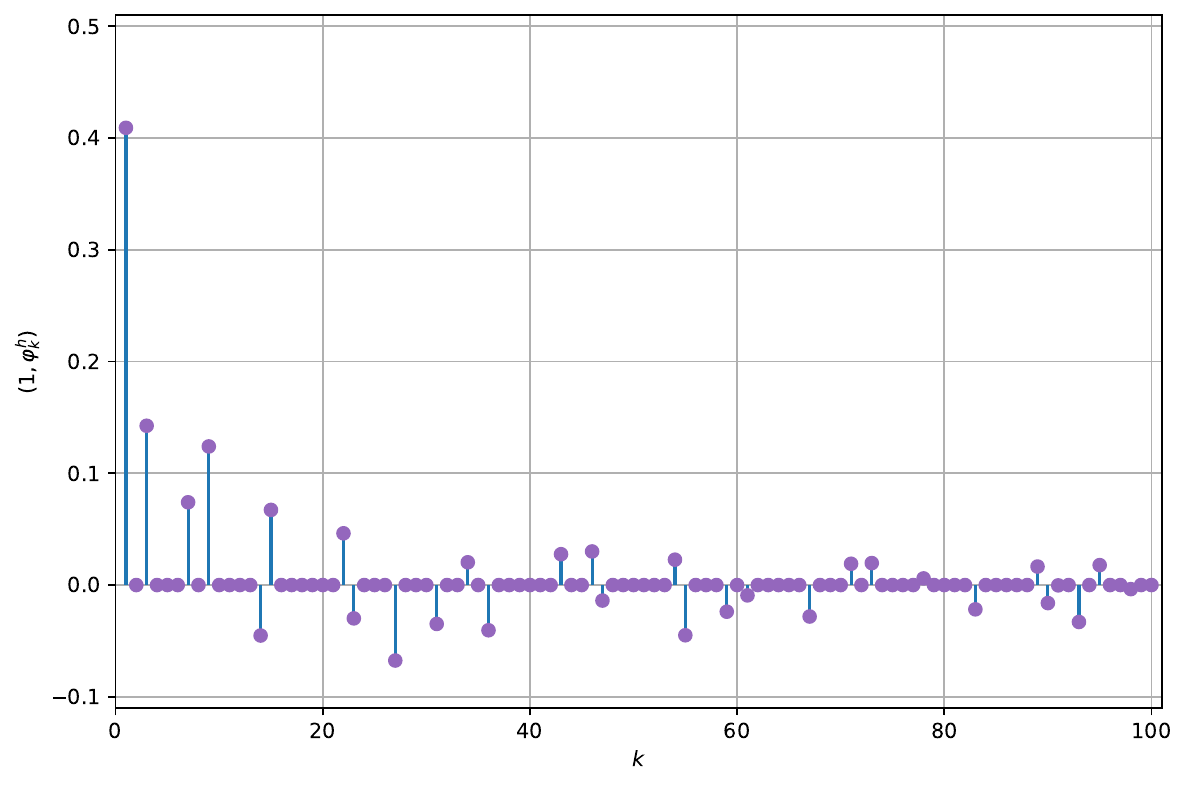}
		\includegraphics[width=0.48\linewidth]{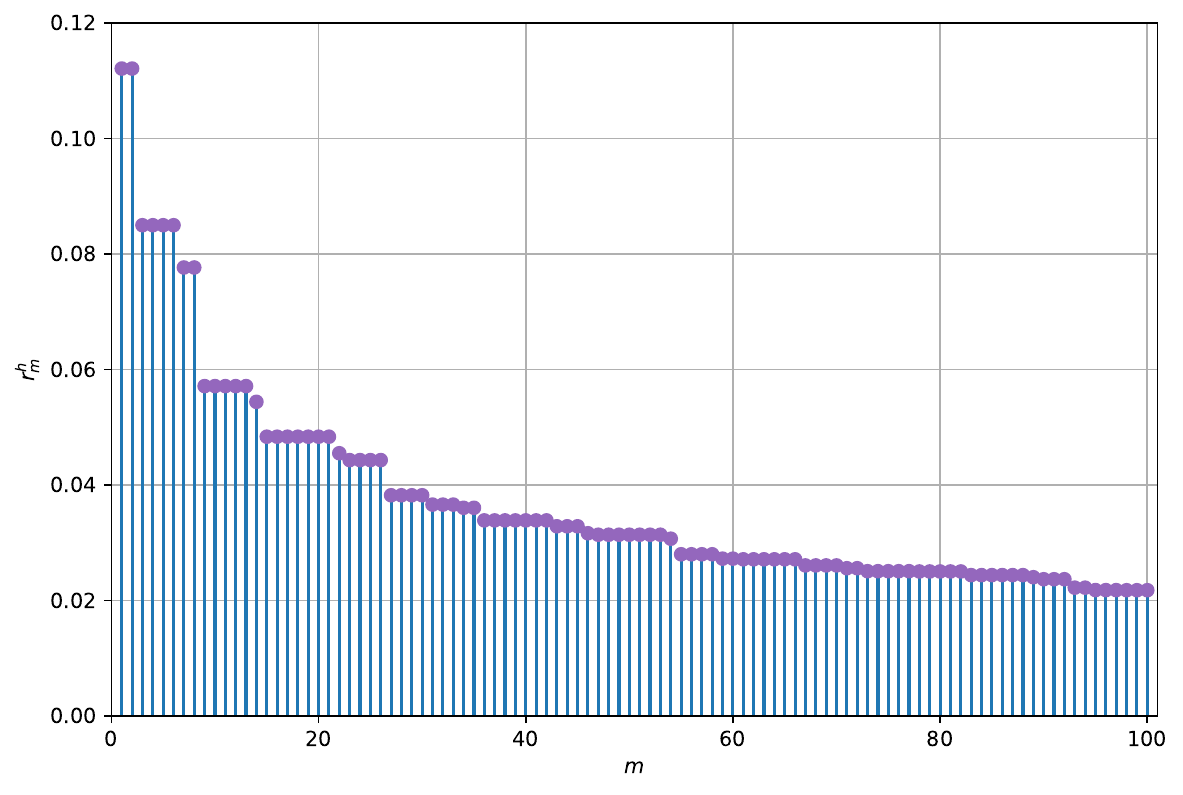}
		\caption{Initial condition expansion coefficients (left) and higher harmonics parameter (right).}
		\label{f-5.8}
	\end{figure}

	When solving the homogenized problem, these approximation errors are suppressed by approximating the kernel contribution from large eigenvalues with a delta-function term. The higher harmonics parameter $r_m^h$ is calculated according to \eqref{3.9}. For our test problem, the results are shown in Fig.~\ref{f-5.8}. Numerical data for selected filtered coefficients of the memory kernel approximation ($m(\epsilon)$, $\epsilon=10^{-5}$) are given in Table~4.
	
	\subsection{Solution of the homogenized problem}
	
	Numerical calculations are first performed for the local model. We obtain this model when the entire memory kernel is approximated by a delta function. In the limiting case $m=0$, from \eqref{3.10} we get the parabolic equation:
	\[
	(1 + r_0) \frac{\partial \widetilde{u}(x,t)}{\partial t} = \nabla \cdot \left(D \nabla \widetilde{u}(x,t)\right), \quad x \in \Omega,\ t > 0.
	\]
	In our case (see \eqref{3.9}):
	\[
	r_0 = (1-|Y_2|)^{-1} |Y_2| \approx 0.335697.
	\]
	
\begin{figure}[ht]
		\centering
		\includegraphics[width=0.48\linewidth]{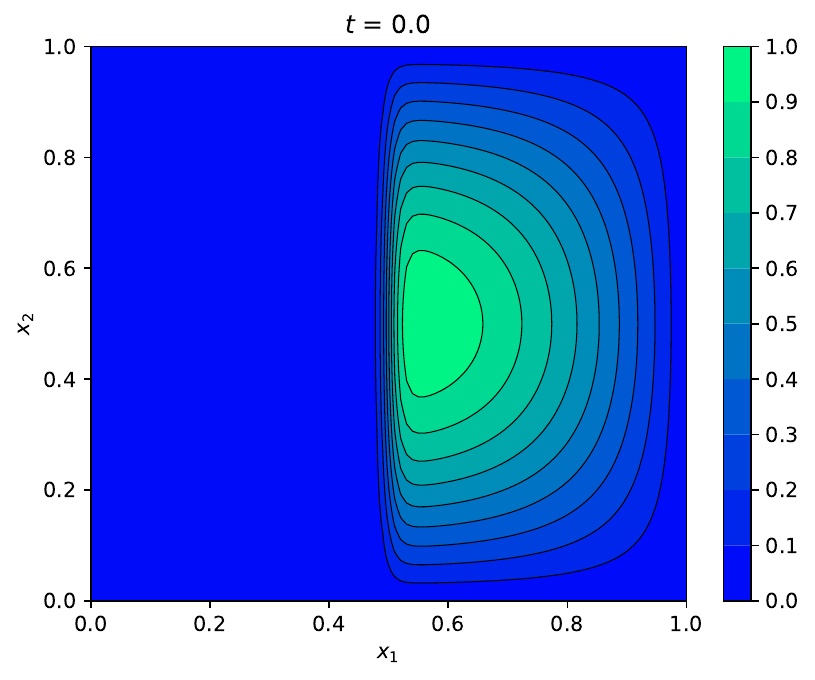}
		\includegraphics[width=0.48\linewidth]{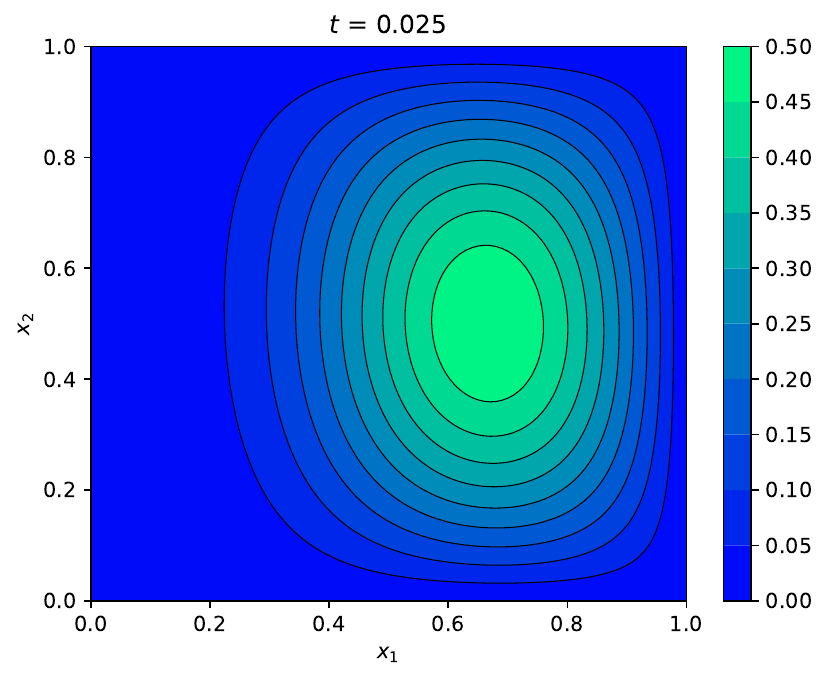} \\
		\includegraphics[width=0.48\linewidth]{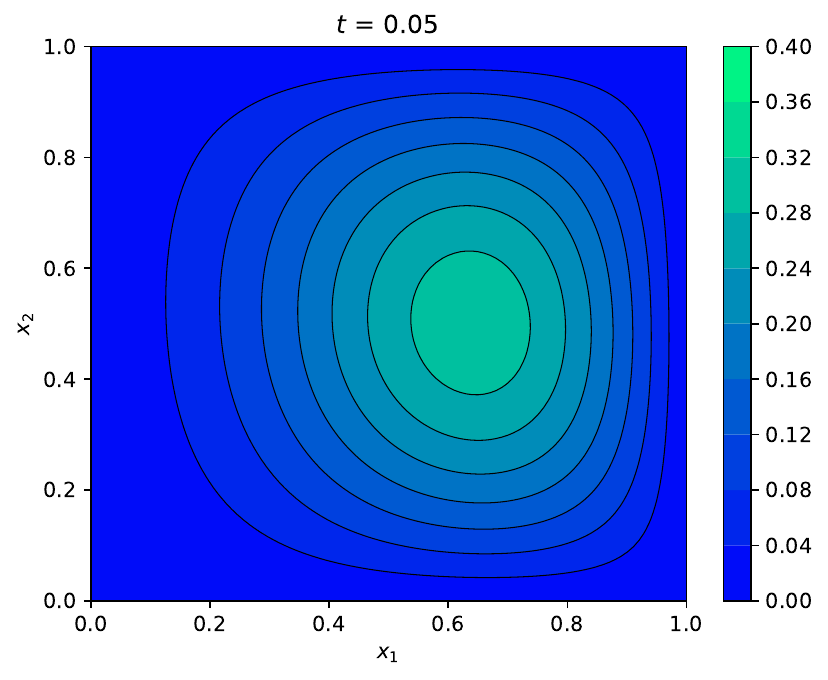}
		\includegraphics[width=0.48\linewidth]{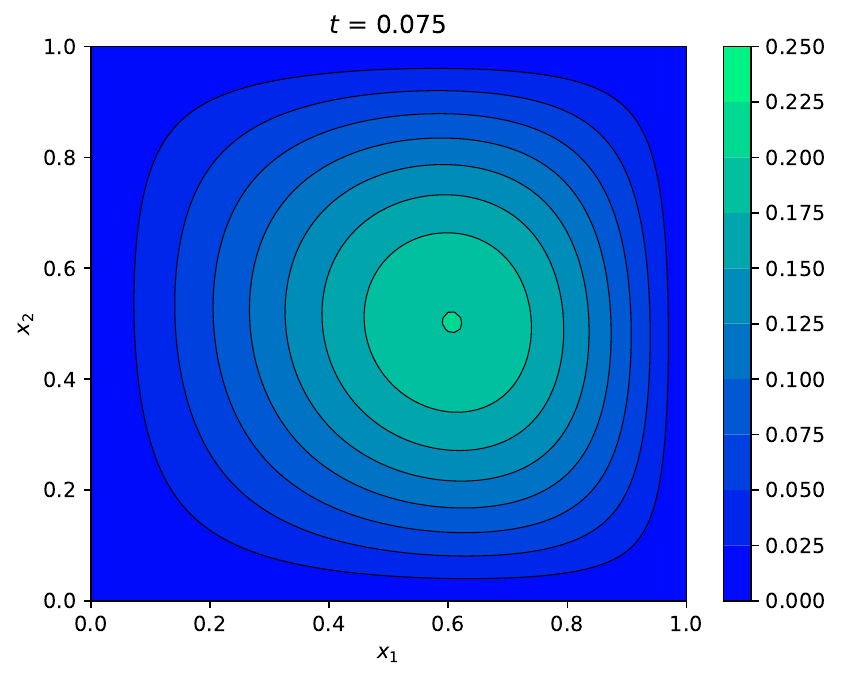} \\
		\caption{Solution of the problem without memory effects at different time moments.}
		\label{f-5.9}
	\end{figure}

\begin{figure}[ht]
		\centering
		\includegraphics[width=0.48\linewidth]{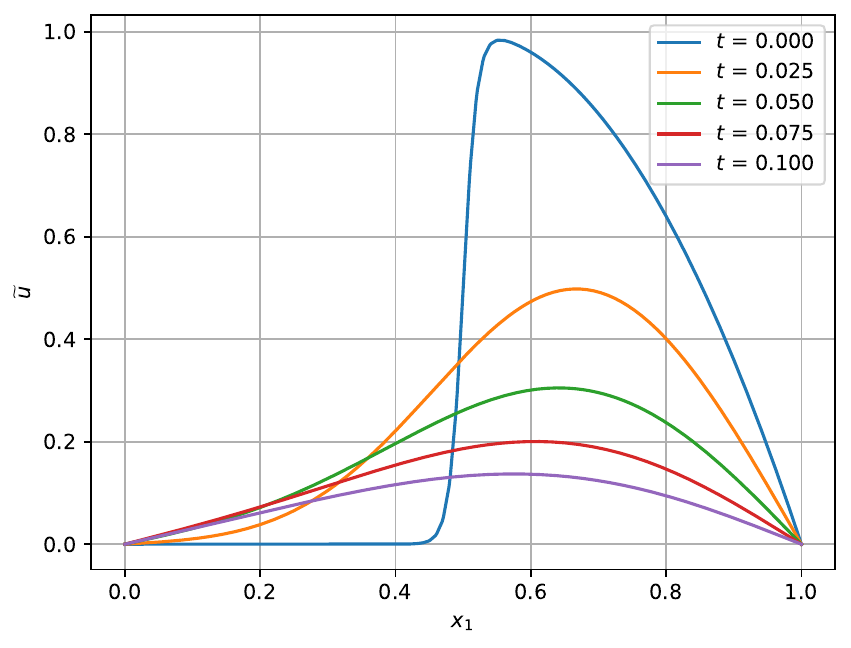}
		\includegraphics[width=0.48\linewidth]{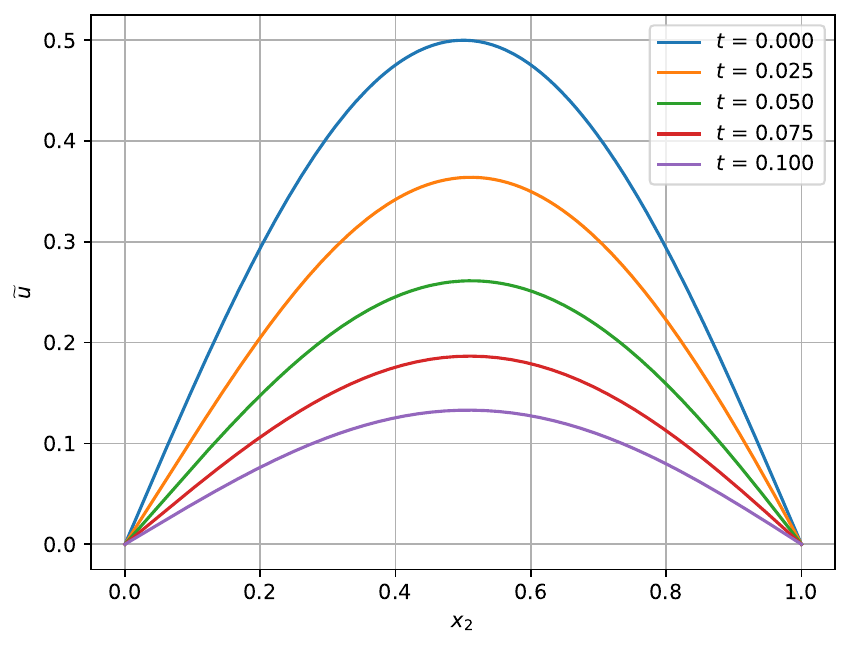} \\
		\caption{Solution of the problem without memory effects in two cross-sections.}
		\label{f-5.10}
	\end{figure}

	In the domain $\Omega$, we construct a triangulation by first partitioning the square uniformly into smaller rectangles of equal size $h \times h$, and then dividing each rectangle into two right triangles using the diagonal from the bottom-left to the top-right vertex.
	The approximate solution $\widetilde{u}^h(x,t)$ at specific time moments is shown in Fig.~\ref{f-5.9}. The solution along the lines $x_1 = 0.5$ and $x_2 = 0.5$ is presented in Fig.~\ref{f-5.10}. The calculations were performed using a fully implicit scheme ($\sigma=1$) with time step $\tau = 10^{-4}$ and uniform spatial discretization with $h=0.01$ in each direction.
	
	\begin{figure}[ht]
		\centering
		\includegraphics[width=0.48\linewidth]{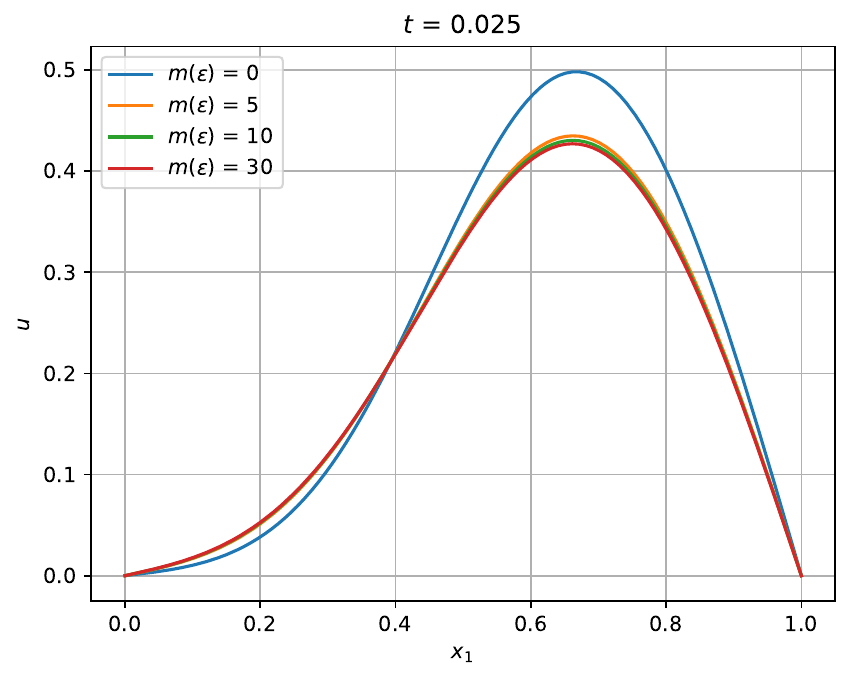}
		\includegraphics[width=0.48\linewidth]{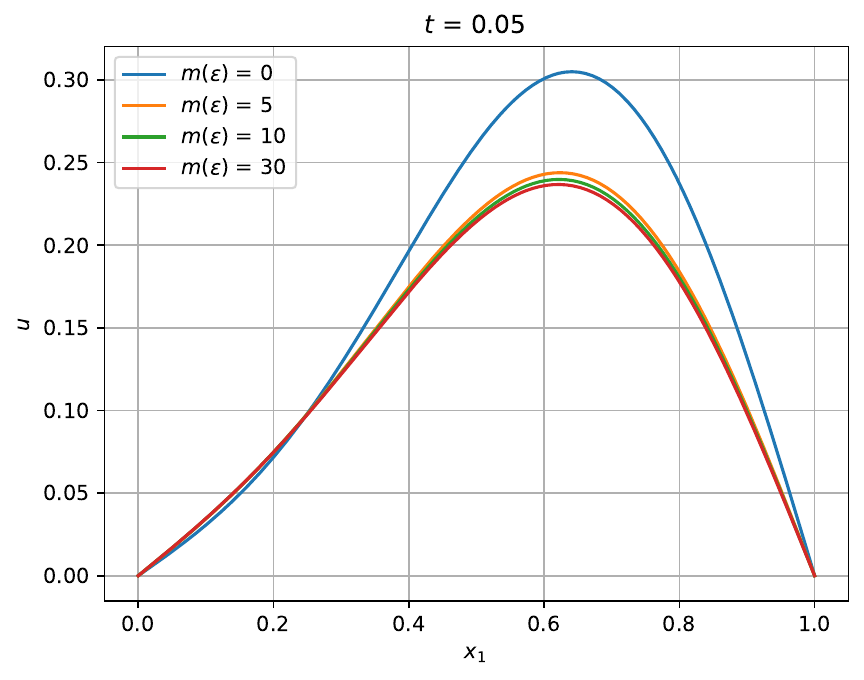} \\
		\includegraphics[width=0.48\linewidth]{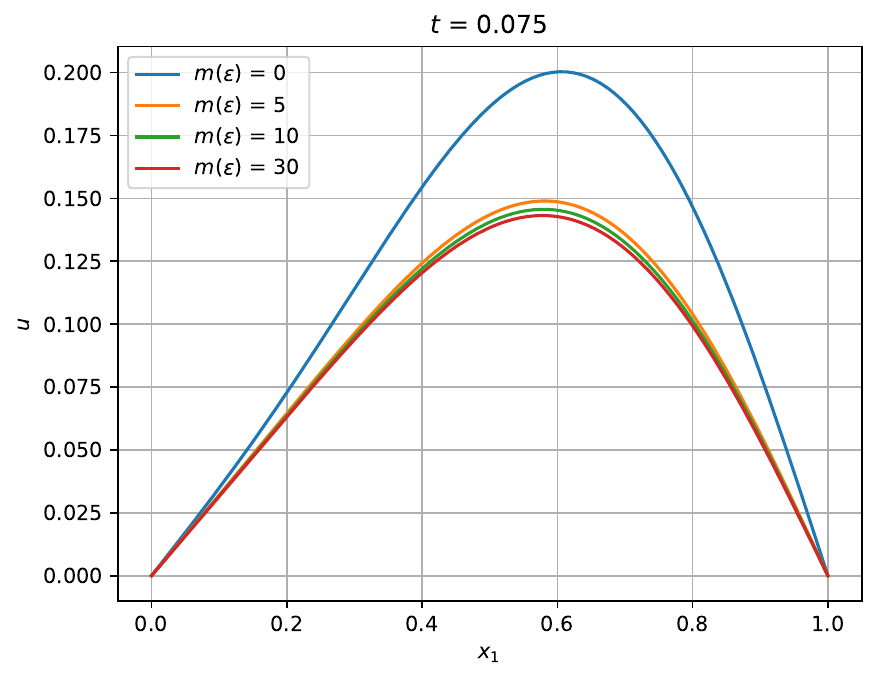}
		\includegraphics[width=0.48\linewidth]{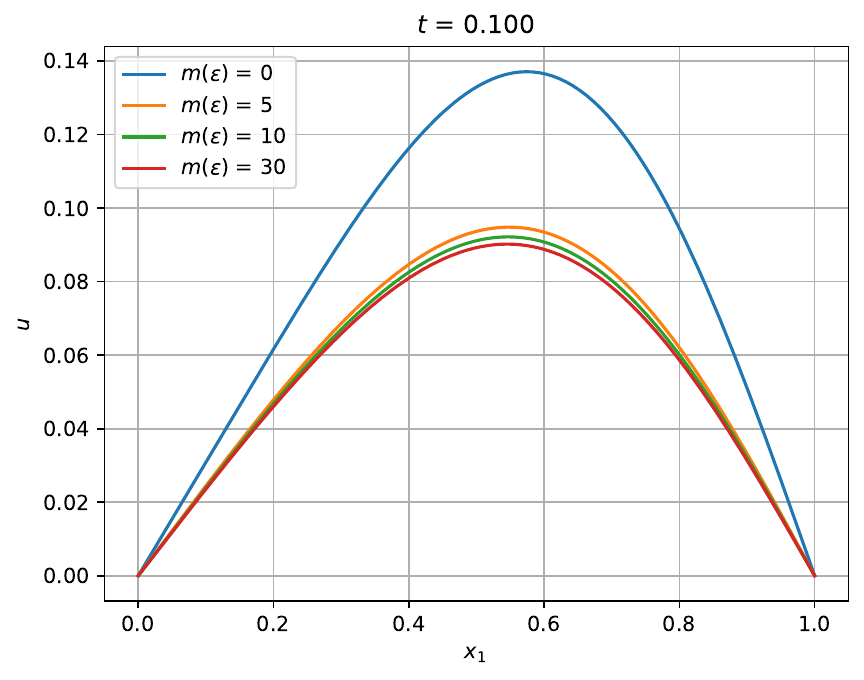} \\
		\caption{Solution of the problem at different time moments.}
		\label{f-5.11}
	\end{figure}
	
	\clearpage
	
	\section{Conclusions} \label{sec:6}
	
	\begin{enumerate}[(1)]
		
		\item We consider a homogenization problem for the nonstationary diffusion equation with weak diffusion in inclusions. The mathematical homogenization model incorporates memory effects described by a Volterra integro-differential equation at the macroscale level. The memory kernel is determined from the solution of a nonstationary problem on the periodicity cell.
		
		\item A novel approach is proposed for solving homogenization problems of nonstationary processes with memory effects. The method is based on: (i) standard determination of the effective diffusion coefficient from solutions to auxiliary stationary boundary-value problems on the periodicity cell; (ii) approximation of the memory kernel by a sum of exponentials obtained through solving a partial spectral problem on the periodicity cell; and (iii) formulation of a local problem at the macroscale level for an extended system of nonstationary equations.
		
		\item The developed computational homogenization algorithm combines: solutions of stationary and nonstationary problems on the periodicity cell with a boundary-value problem for a system of nonstationary equations at the macroscale level. Spatial approximation employs Lagrange finite elements, while temporal discretization uses two-level weighted schemes. Stability estimates for the approximate solution in appropriate Hilbert spaces are established, consistent with estimates for the differential problem solution.
		
		\item Key elements of the computational algorithm and its effectiveness are demonstrated through numerical solution of a two-dimensional homogenization problem. Special attention is given to the problem of approximating the memory kernel by exponential sums.
		
	\end{enumerate}


\begin{thebibliography}{10}
	\expandafter\ifx\csname url\endcsname\relax
	\def\url#1{\texttt{#1}}\fi
	\expandafter\ifx\csname urlprefix\endcsname\relax\def\urlprefix{URL }\fi
	\expandafter\ifx\csname href\endcsname\relax
	\def\href#1#2{#2} \def\path#1{#1}\fi
	
	\bibitem{efendiev2009multiscale}
	Y.~Efendiev, T.~Y. Hou, Multiscale Finite Element Methods: Theory and
	Applications, Springer Science \& Business Media, 2009.
	
	\bibitem{pavliotis2008multiscale}
	G.~A. Pavliotis, A.~Stuart, Multiscale Methods: Averaging and Homogenization,
	Springer Science \& Business Media, 2008.
	
	\bibitem{steinhauser2017computational}
	M.~Steinhauser, Computational Multiscale Modeling of Fluids and Solids,
	Springer, 2017.
	
	\bibitem{bakhvalov2012homogenisation}
	N.~S. Bakhvalov, G.~Panasenko, Homogenisation: Averaging Processes in Periodic
	Media: Mathematical Problems in the Mechanics of Composite Materials,
	Springer Science \& Business Media, 2012.
	
	\bibitem{bensoussan2011asymptotic}
	A.~Bensoussan, J.-L. Lions, G.~Papanicolaou, AsYmptotic Analysis for Periodic
	Structures, American Mathematical Soc., 2011.
	
	\bibitem{sandrakov1}
	G.~V. Sandrakov, Homogenization of parabolic equations with contrasting
	coefficients, Izv. Math. 63~(5) (1999) 1015–1061.
	
	\bibitem{sandrakov2}
	G.~V. Sandrakov, V.~V. Semenov, Homogenized models with memory effect for
	heterogeneous periodic media, Journal of Optimization, Differential Equations
	and Their Applications 30~(2) (2022) 1--18.
	
	\bibitem{ChenBook1998}
	C.~Chen, T.~Shih, Finite Element Methods for Integrodifferential Equations,
	World Scientific, Singapore, 1998.
	
	\bibitem{linz1985analytical}
	P.~Linz, Analytical and Numerical Methods for Volterra Equations, SIAM, 1985.
	
	\bibitem{vabishchevich2022numerical}
	P.~N. Vabishchevich, Numerical solution of the {Cauchy} problem for {Volterra}
	integrodifferential equations with difference kernels, Applied Numerical
	Mathematics 174 (2022) 177--190.
	
	\bibitem{vabishchevich2023approximate}
	P.~N. Vabishchevich, Approximate solution of the {Cauchy} problem for a
	first-order integrodifferential equation with solution derivative memory,
	Journal of Computational and Applied Mathematics 422 (2023) 114887.
	
	\bibitem{braess1986nonlinear}
	D.~Braess, Nonlinear Approximation Theory, Springer, Berlin, Heidelberg, 1986.
	
	\bibitem{beylkin2005approximation}
	G.~Beylkin, L.~Monz{\'o}n, On approximation of functions by exponential sums,
	Applied and Computational Harmonic Analysis 19~(1) (2005) 17--48.
	
	\bibitem{hornung1997homogenization}
	U.~Hornung (Ed.), Homogenization and Porous Media, Springer Science \& Business
	Media, 1997.
	
	\bibitem{KnabnerAngermann2003}
	P.~Knabner, L.~Angermann, Numerical Methods for Elliptic and Parabolic Partial
	Differential Equations, Springer, New York, 2003.
	
	\bibitem{QuarteroniValli1994}
	A.~Quarteroni, A.~Valli, Numerical Approximation of Partial Differential
	Equations, Springer-Verlag, Berlin, 1994.
	
	\bibitem{Ascher2008}
	U.~M. Ascher, Numerical Methods for Evolutionary Differential Equations,
	Society for Industrial and Applied Mathematics, 2008.
	
	\bibitem{SamarskiiTheory}
	A.~A. Samarskii, The Theory of Difference Schemes, Marcel Dekker, New York,
	2001.
	
	\bibitem{mclean1993numerical}
	W.~McLean, V.~Thom{\'e}e, Numerical solution of an evolution equation with a
	positive-type memory term, The ANZIAM Journal 35~(1) (1993) 23--70.
	
	\bibitem{mclean1996discretization}
	W.~McLean, V.~Thom{\'e}e, L.~B. Wahlbin, Discretization with variable time
	steps of an evolution equation with a positive-type memory term, Journal of
	Computational and Applied Mathematics 69~(1) (1996) 49--69.
	
	\bibitem{evans2010partial}
	L.~C. Evans, Partial Differential Equations, American Mathematical Society,
	2010.
	
	\bibitem{halanay1965asymptotic}
	A.~Halanay, On the asymptotic behavior of the solutions of an
	integro-differential equation, Journal of Mathematical Analysis and
	Applications 10~(2) (1965) 319--324.
	
	\bibitem{saad2011numerical}
	Y.~Saad, Numerical Methods for Large Eigenvalue Problems: Revised Edition,
	SIAM, 2011.
	
	\bibitem{SamarskiiMatusVabischevich2002}
	A.~A. Samarskii, P.~P. Matus, P.~N. Vabishchevich, Difference Schemes with
	Operator Factors, Kluwer Academic, Dordrecht, 2002.
	
	\bibitem{geuzaine2009gmsh}
	C.~Geuzaine, J.-F. Remacle, Gmsh: A {3-D} finite element mesh generator with
	built-in pre-and post-processing facilities, International journal for
	numerical methods in engineering 79~(11) (2009) 1309--1331.
	
	\bibitem{alnaes2015fenics}
	M.~Aln{\ae}s, J.~Blechta, J.~Hake, A.~Johansson, B.~Kehlet, A.~Logg,
	C.~Richardson, J.~Ring, M.~E. Rognes, G.~N. Wells, The {FEniCS} project
	version 1.5, Archive of numerical software 3~(100).
	
\end{thebibliography}
	\end{document}